\documentclass[12pt]{article}
\usepackage{amsmath}
\usepackage{amssymb}

\newtheorem{thm}{Theorem}[section] 
\newtheorem{lemma}[thm]{Lemma} 
\newtheorem{Def}[thm]{Definition} 
 
\newtheorem{prop}[thm]{Proposition}

\newtheorem{examp}[thm]{Example}
\newtheorem{hyp}[thm]{Hypothesis}

\newcommand{\proof}{\noindent {\bf Proof} \hspace{0.2in}} 
\newcommand{\qed}{\hfill\mbox{\raggedright\rule{.07in}{.1in}}
  \vspace{1ex}} 
\newcommand{\dps}{\displaystyle}
\newcommand{\Section}[1]{\section{#1} \setcounter{equation}{0}}

\title{Toroidal normal forms for bifurcations in retarded functional
differential equations II: Saddle-node/multiple Hopf interaction}
 
\author{Younsun Choi\\Department of Mathematics and Statistics\\
University of Ottawa\\Ottawa, ON K1N 6N5\\CANADA  
\and Victor G. LeBlanc\\Department of
  Mathematics and Statistics\\University of Ottawa\\Ottawa, 
ON K1N 6N5\\CANADA}
\date{July 4, 2005}

\begin{document}

\maketitle

\begin{abstract}
In this paper, we study the realizability problem for retarded
 functional
differential equations near an equilibrium point undergoing
a nonlinear mode interaction between a
saddle-node bifurcation and a non-resonant multiple Hopf bifurcation.  
In contrast to the
case of transcritical/multiple Hopf interaction which was studied in
an earlier paper \cite{CL1}, the analysis here is complicated by the presence of
a nilpotency which introduces a non-compact component in the
symmetry group of the normal form.  
We present a framework to analyse
the realizability problem in this non-semisimple case which exploits
to a large extent our previous results for the realizability problem
in the semisimple case.  
Apart from providing a solution to the problem of interest in this paper,
it is believed that the approach used here could potentially be
 adapted to the
study of the realizability problem for toroidal normal forms in the general case of repeated
eigenvalues with Jordan blocks.
\end{abstract}

\pagebreak
\Section{Introduction}
Retarded functional differential equations (RFDEs) are frequently used as models for various
phenomena \cite{BC94,BBL,HFEKGG,Kuang,LM,SC,SS}.  
While the phase space for the resulting dynamical
system is infinite dimensional, the existence of finite-dimensional
invariant center manifolds near bifurcation points imply that much of
the machinery developed for the analysis of finite codimension
bifurcations in ordinary differential equations (e.g. normal forms, unfolding theory) are portable to
RFDEs.  Indeed, there is ample evidence in the literature, e.g. 
\cite{BC94,HFEKGG,SC}, that
these tools and techniques of local analysis can give much valuable information
about the dynamics of RFDEs.

In this context, one of the fundamental questions concerns the
characterization of the range of dynamics accessible near a bifurcation
point for a given
RFDE model.  This is not a trivial question, since even for scalar RFDEs, 
it is the possible to have bifurcations of equilibria with
large dimensional center manifolds.
In this case, if there are not enough independent delay terms in the
nonlinear part of the RFDE, there may be severe restrictions on the
possible dynamics which can be realized in the center manifold
equations.  The study of these types of questions is known generally as the {\em
  realizability problem} for RFDEs, and we refer the reader to the
Introduction of \cite{CL1} and to \cite{FMR,FM96,Hal85} for more details.

In a previous paper \cite{CL1}, we studied the realizability problem for scalar
RFDEs in two important cases: the case of multiple non-resonant Hopf
bifurcation, and the case of the interaction between a transcritical
steady-state bifurcation and a multiple non-resonant Hopf
bifurcation.  In particular, we used the fact that these bifurcations
admit normal forms with toroidal symmetries, and that these
symmetries allow for a decoupling of the center manifold normal form equations
into a {\em radial part} and an {\em angular part}.  The uncoupled
radial equations are a crucial part of the dynamics near the
bifurcation, and it is thus very reasonable to consider the
problem of realizability of the class of radial equations within a
given RFDE.  This problem was solved in general in \cite{CL1} for the two
important bifurcation scenarios described above.
Moreover, it was shown in \cite{CL1} that our solution to this
realizability
problem, that is, our estimate on the
number of independent delays sufficient to achieve complete realizability, is optimal.
The results in \cite{CL1}
considerably generalize previous results of Faria and
Magalh$\tilde{\mbox{\rm a}}$es \cite{FM96} and 
of Buono and B\'elair \cite{BB}.

However, in \cite{CL1} 
the case of the interaction between a saddle-node
steady-state bifurcation and a multiple non-resonant Hopf bifurcation
was not studied, because there is a nilpotency associated to this
bifurcation which is such that the approach used in 
\cite{CL1} is not applicable.
This nilpotency arises from the fact that the normal form for the saddle-node
bifurcation contains an unfolding parameter which appears as an affine
linear perturbation of the singular vector field, i.e. 
\[
\begin{array}{lll}
\dot{\rho}&=&\nu+a\rho^2+O(\rho^3)\\[0.1in]
\dot{\nu}&=&0.
\end{array}
\]
A consequence of this nilpotency is that the normal form 
for the saddle-node/multiple Hopf interaction
admits a
symmetry group which is a direct product of a compact torus group, and
a non-compact group isomorphic to the additive group of real numbers.
This latter non-compact portion does not occur in the cases studied in 
\cite{CL1}.
It thus becomes important to develop a framework which allows for the
characterization of the effects of this
non-compact symmetry on the normal form before being able to address
the realizability problem for this bifurcation in RFDEs.  

In this paper, we adapt the Faria and Magalh$\tilde{\mbox{\rm a}}$es
normal form procedure \cite{FMTB,FMH} to the saddle-node/multiple Hopf
interaction in scalar RFDEs, and use it to develop a framework suitable to
studying the realizability problem for the radial part of the
normal form in this nilpotent (non-semisimple) case.  The framework is
carefully constructed so that in the end, we may use our solution to the
realizability problem in the semisimple case in \cite{CL1} to the fullest
possible extent.

\Section{Preliminaries}
In this section we will briefly recall some standard results and
terminology in the bifurcation theory of RFDEs in order
to establish the notation.  For more details, see \cite{CL1,FMTB,FMH,HalVL}.

\subsection{Infinite dimensional parameterized ODE}
Suppose $r>0$ is a given real
number, $n\geq 1$ is a given integer and $C_n\equiv C\left(  \left[  -r,0\right]  ,\mathbb{R}^n\right)  $ is the
Banach space of continuous functions from $\left[  -r,0\right]  $ into
$\mathbb{R}^n$ with supremum norm.  We define $u_{t}\in C_n$ as $u_{t}\left(  \theta\right)
=u\left(  t+\theta\right)  ,-r\leq\theta\leq 0.$ 
Let us consider the
following parameterized family of scalar $(n=1)$ nonlinear retarded functional differential equations
\begin{equation}
\dot{z}\left(  t\right)  =L(\alpha)z_{t}+F\left(  z_{t},\alpha\right)  ,
\label{y1}%
\end{equation}
\noindent where $L:C_1\times\mathbb{R}^{s+1}\rightarrow\mathbb{R}$ is
a parameterized family ($s\geq 0$) of bounded linear
operators from $C_1$ into $\mathbb{R}$ and $F$ is
a smooth function from $C_{1}\times\mathbb{R}^{s+1}$ into
$\mathbb{R}$.
In the ``prequel'' \cite{CL1} to this paper, we assumed that $F(0,0)=0$, and
$DF(0,0)=0$.  While this hypothesis included the cases of multiple
Hopf bifurcation and of a mode-interaction between a multiple Hopf and
a transcritical type steady-state bifurcation, 
it excluded the case of a mode interaction in which the steady-state was of
saddle-node type.  Therefore, in this paper, we assume the
following weaker hypothesis
\begin{hyp}
$F(0,0)=0,\,\,D_1F(0,0)=0,\,\,\,\mbox{\rm and}\,\,D_{\alpha}F(0,0)\neq
0$.
\label{hyp1}
\end{hyp}
Performing a linear change of parameters and relabeling the parameters
if necessary, we may then rewrite (\ref{y1}) as
\begin{equation}
\dot{z}(t)=L_0z_t+\nu+\hat{F}(z_t,\nu,\mu),
\label{y1split}
\end{equation}
where we have set $\alpha\equiv
(\nu,\mu)\in\mathbb{R}^{1}\times\mathbb{R}^s$, $L_0\equiv L(0)$, and
$\hat{F}(z_t,\nu,\mu)=(L(\nu,\mu)-L_0)z_t+F(z_t,\nu,\mu)-D_{\nu}F(0,0,0)$.  
It follows from Hypothesis \ref{hyp1} that
\[
\hat{F}(0,0,0)=0,\,\,\,\,\,\,\,\,\,D\hat{F}(0,0,0)=0.
\]
Clearly, the
parameter $\nu$ plays a distinguished role in (\ref{y1split}) in
comparison to the other $s$ parameters $\mu$.  

\vspace*{0.15in}
\noindent
{\em Spectral hypothesis}
\vspace*{0.15in}

Suppose we set $(\nu,\mu)=(0,0)$ in (\ref{y1split}).  
If we then write $L_0$ as
\begin{equation}
L_0\phi=\int_{-r}^{0}\,d\eta\left(\theta\right)\phi\left(
\theta\right)  ,
\label{L0def}
\end{equation}
where $\eta$ is a real-valued function of bounded variation in
$[-r,0]$ and we let $A_0$ be the infinitesimal generator of the
semi-flow associated with the linear RFDE $\dot{z}(t)=L_0z_t$, then it
is well-known that the spectrum $\sigma(A_0)$ of $A_0$ is equal to the
point spectrum of $A_0$, and $\lambda\in\sigma(A_0)$ if and only if
$\lambda$ satisfies the characteristic equation
\begin{equation}
\det\Delta\left(  \lambda\right)  =0,\text{ \ \ \ }\Delta\left(  \lambda\right)
=\lambda-\int\limits_{-r}^{0}\,d\eta\left(\theta\right)
e^{\lambda\theta}. \label{ychar}%
\end{equation}
Denote by
$\Lambda_0$ the set of eigenvalues of $\sigma(A_0)$ with zero
real part.
\begin{hyp}
Throughout the rest of the paper, we assume the following hypotheses on $\Lambda_0$.
Each element of $\Lambda_0$ is a simple eigenvalue of
  $A_0$, and $\Lambda_0$ has the following form:
\[
\Lambda_0=\{0,\pm\,i\omega_1,\ldots,\pm\,i\omega_p\},
\]
where
$\omega_1,\ldots,\omega_p$, are independent over the rationals,
i.e. if $r_1,\ldots,r_p$ are rational numbers such that
$
\sum_{j=1}^p\,r_j\omega_j=0
$, then $r_1=\cdots=r_p=0$.  We further assume that the rest of the
spectrum of $A_0$ is bounded away from the imaginary axis.
\label{spectralhyp}
\end{hyp}

\vspace*{0.15in}
\noindent
{\em Phase space decomposition}
\vspace*{0.15in}

In order to properly analyse the role of the parameters in equation
(\ref{y1split}) we need to augment this equation by considering the following
system
\begin{subequations}\label{y1sys:gp}
\begin{align}
\left(\begin{array}{c}\dot{z}(t)\\\dot{\nu}(t)\end{array}\right)&=
\left(\begin{array}{c}L_0z_t+\nu(0)+{\cal F}\left(\left(\begin{array}{c}z_t\\\nu_t\end{array}\right),\mu_t\right)\\0\end{array}\right)\label{y1sys:gp1}\\
\dot{\mu}(t)&=0,\label{y1sys:gp2}
\end{align}
\end{subequations}
where ${\cal F}((z_t,\nu_t)^T,\mu)=\hat{F}(z_t,\nu(0),\mu(0))$.
The above system (\ref{y1sys:gp}) can thus be viewed as an $s$-dimensional family
(parameterized by $\mu$) of {\em 2-dimensional} RFDEs.

Taking into account Hypothesis \ref{spectralhyp}, the linearized
equation $(\dot{z}(t),\dot{\nu}(t))=(L_0z_t+\nu(0),0)$ associated to
(\ref{y1sys:gp1})
has simple non-resonant characteristic values $\pm\,
i\omega_1,\ldots,\pm\, i\omega_p$, and a characteristic value at $0$
of multiplicity $2$. 
We then let $P\subset C_2$ designate the $2p+2$-dimensional center
subspace which is spanned by the columns of the following matrix
\begin{equation}
\mathring{\Phi}=\left(\begin{array}{ccccccc}
1&e^{i\omega_1\theta}&e^{-i\omega_1\theta}&\cdots&e^{i\omega_p\theta}&e^{-i\omega_p\theta}&\theta\\
0&0&0&\cdots&0&0&1
\end{array}\right).
\label{phichoice}
\end{equation}
We note that $\mathring{\Phi}$ satisfies the linear differential equation
${\dps\frac{d\mathring{\Phi}}{d\theta}=\mathring{\Phi}\mathring{B}}$,
where $\mathring{B}$ is the $(2p+2)\times (2p+2)$ matrix
$\mathring{B}=\mbox{\rm diag}(0,i\omega_1,-i\omega_1,\ldots,i\omega_p,-i\omega_p,0)+\mathring{N}$,
and $\mathring{N}$ is the nilpotent matrix whose entries are all zero except the
entry at the intersection of the first row and last column, whose
value is 1.

We decompose $C_2$ as
\begin{equation}
C_2=P\oplus Q.
\label{C1split}
\end{equation}
Defining $C^*_2\equiv
C([0,r],\mathbb{R}^{2*})$, where $\mathbb{R}^{2*}$ is the
2-dimensional space of row vectors, we introduce the adjoint bilinear form on
$C_2^*\times C_2$:
\begin{equation}
(\psi,\phi)=\psi(0)\phi(0)-\int_{-r}^0\,\int_0^{\theta}\,\psi(\xi-\theta)d\tilde{\eta}(\theta)\phi(\xi)d\xi,
\label{bf}
\end{equation}
where $d\tilde{\eta}(\theta)$ is of the form
\[
d\tilde{\eta}(\theta)=\left(\begin{array}{cc}d\eta(\theta)&\delta(\theta)\,d\theta\\0&0\end{array}\right),
\]
$d\eta(\theta)$ is as in (\ref{L0def}) and $\delta(\theta)\,d\theta$ is
such that
\[
\int_{-r}^0\,f(\theta)\,\delta(\theta)\,d\theta=f(0).
\]

We may then choose a basis $\{\psi_1,\ldots,\psi_{2p+2}\}$ of the dual
space $P^*$ such
that if $\mathring{\Psi}=\mbox{\rm col}(\psi_1,\ldots,\psi_{2p+2})$ then
$(\mathring{\Psi},\mathring{\Phi})=I_{2p+2}$, where 
throughout the paper we use the convention that for a given integer
$q\geq 1$, $I_{q}$ is the $q\times q$
identity matrix.
The following result will be useful later.
\begin{lemma}
Let $\mathring{\Phi}$ be as in (\ref{phichoice}) and $\mathring{\Psi}(\xi)=(\psi_{k,\ell}(\xi))$ be a
$(2p+2)\times 2$ matrix whose rows form a basis for $P^*$, and such
that $(\mathring{\Psi},\mathring{\Phi})=I_{2p+2}$.
Then
\[
\psi_{2p+2,1}(0)=0,\,\,\,\,\,\mbox{\rm and}\,\,\,\,\,\psi_{k,1}(0)\neq
0,\,1\leq k\leq 2p+1.
\]
\label{psi0lem}
\end{lemma}
\proof
It follows from Hypothesis \ref{spectralhyp} that 
\[
\int_{-r}^0\,d\eta(\theta)=0,\,\,\,\,\int_{-r}^0\,d\eta(\theta)e^{\pm\,i\omega_j\theta}=\pm\,i\omega_j,\,\,\,j=1,\ldots,p
\]
and
\[
\int_{-r}^0\,d\eta(\theta)\,\theta\neq
1,\,\,\,\,\int_{-r}^0\,d\eta(\theta)\,\theta\,e^{\pm\,i\omega_j\theta}\neq
1,\,\,\,j=1,\ldots,p.
\]
Using (\ref{bf}), the equation $(\Psi,\Phi)=I_{2p+2}$ leads to
\[
\psi_{2p+2,1}(0)\left(1-\int_{-r}^0\,d\eta(\theta)\,\theta\right)=0,
\]
\[
\psi_{k,1}(0)\left(1-\int_{-r}^0\,d\eta(\theta)\,\theta\,e^{i\omega_j\theta}\right)=1,\,\,\,k=2,4,\ldots,2p,
\]
and
\[
\psi_{k,1}(0)\left(1-\int_{-r}^0\,d\eta(\theta)\,\theta\,e^{-i\omega_j\theta}\right)=1,\,\,\,k=3,5,\ldots,2p+1,
\]
so we get the desired conclusion.
\hfill\qed

Consider now the Banach space $BC_2$ of functions from $[-r,0]$
into $\mathbb{R}^2$ which are uniformly continuous on $[-r,0)$ with a
jump discontinuity at 0.  Elements of $BC_2$ are written as
\[
\phi+X_0\lambda
\]
where $\phi\in C_2$, $\lambda\in\mathbb{R}^2$ and $X_0$ is the $2\times
2$ matrix-valued function
\[
X_0(\theta)=\left\{\begin{array}{cc}
I_2&\theta=0\\&\\
0&-r\leq\theta <0.
\end{array}\right.
\]
Let $\pi: BC_2\rightarrow P$ denote the projection
\[
\pi(\phi+X_0\lambda)=\mathring{\Phi}\,[(\mathring{\Psi},\phi)+\mathring{\Psi}(0)\lambda],
\]
where $(\,\,\,,\,\,\,)$
is the bilinear form (\ref{bf}).  We may then extend the splitting
(\ref{C1split}) to
\begin{equation}
BC_2=P\oplus\mbox{\rm ker}\,\pi,
\label{BC1split}
\end{equation}
with the property that $Q\subsetneq\mbox{\rm ker}\,\pi$.

This structure now allows for a decomposition of the
phase space which facilitates
the implementation of the normal form procedure
with parameters for RFDEs developed by Faria and
Magalh$\tilde{\mbox{\rm a}}$es in \cite{FMH}.
Specifically, it follows that (\ref{y1sys:gp1}) is equivalent to the
parameterized family
\begin{subequations}\label{y3}
\begin{align}
{\dps\left(\begin{array}{c}\dot{x}\\\dot{\nu}\end{array}\right)}&=
{\dps
  \mathring{B}\left(\begin{array}{c}x\\\nu\end{array}\right)+\mathring{\Psi}(0)\left(\begin{array}{c}
{\cal F}\left(\mathring{\Phi}\,\left(\begin{array}{c}x\\\nu\end{array}\right)+y,\mu\right)\\0\end{array}\right)}\label{y3:1}\\&\notag\\
{\dps\frac{d}{dt}y}&={\dps A_{Q^1}y+(I-\pi)X_0\left(\begin{array}{c}
{\cal F}\left(\mathring{\Phi}\,\left(\begin{array}{c}x\\\nu\end{array}\right)+y,\mu\right)\\0\end{array}\right)}\label{y3:2}
\end{align}
\end{subequations}
where $\mathring{\Psi}(0)$ is as in Lemma \ref{psi0lem},
$x\in\mathbb{R}^{2p+1}$, $\nu\in\mathbb{R}$, $\mu\in\mathbb{R}^s$,
$y\in Q^1\equiv
Q\cap C^1_2$, ($C^1_2$ is the subspace of $C_2$ consisting of
continuously differentiable functions), and $A_{Q^1}$ is the operator
from $Q^1$ into $\mbox{\rm ker}\,\pi$ defined by
\[
A_{Q^1}\phi=\dot{\phi}+X_0\,\left[\int_{-r}^0\,d\tilde{\eta}(\theta)\,\phi(\theta)-\dot{\phi}(0)\right].
\]

\subsection{Faria and Magalh$\tilde{\mbox{\bf a}}$es normal form}

Consider the formal Taylor expansion of the nonlinearity ${\cal F}$ terms in (\ref{y1sys:gp})
\[
{\cal F}(u)=\sum_{j\geq 2}\,\frac{1}{j!}\widehat{{\cal F}}_j(u),\,\,\,\,\,u\in\,C_{2+s},
\]
where $\widehat{{\cal F}}_j(w)=H_j(w,\ldots,w)$, with $H_j$ belonging to the space of
continuous multilinear symmetric maps from
$C_{2+s}\times\cdots\times
C_{2+s}$
($j$ times) to $\mathbb{R}$. 
If we denote $\mathbf{x}=(x,\nu)^T$ and $f_j=(f_j^1,f_j^2)$, where
\begin{equation}
\begin{array}{rcl}
f_{j}^1(\mathbf{x},y,\mu)&=&\mathring{\Psi}(0)\,\left(\begin{array}{c}\widehat{{\cal
    F}}_j\left(\mathring{\Phi}\,\mathbf{x}+y,\mu\right)\\0\end{array}\right)\\[0.15in]
f_j^2(\mathbf{x},y,\mu)&=&(I-\pi)\,X_0\,\left(\begin{array}{c}\widehat{{\cal
    F}}_j(\mathring{\Phi}\,\mathbf{x}+y,\mu)\\0\end{array}\right),
\end{array}
\label{fjdef}
\end{equation}
then (\ref{y3}) can be written as
\begin{subequations}\label{y4}
\begin{align}
\dot{\mathbf{x}}&={\displaystyle \mathring{B}\mathbf{x}+\sum_{j\geq
    2}\,\frac{1}{j!}\,f_j^1(\mathbf{x},y,\mu)}\label{y4:1}\\[0.15in]
{\displaystyle\frac{d}{dt}\,y}&={\displaystyle A_{Q^1}y+\sum_{j\geq
    2}\,\frac{1}{j!}\,f_j^2(\mathbf{x},y,\mu)}\label{y4:2}
\end{align}
\end{subequations}

The spectral hypotheses we have specified in Hypothesis
\ref{spectralhyp} are sufficient to conclude that the non-resonance
condition of Faria and Magalh$\tilde{\mbox{\rm a}}$es \cite{FMH} 
holds.  Consequently, using
successively at each order $j$ a near identity change of variables of the form
\begin{equation}
(\mathbf{x},y)=(\hat{\mathbf{x}},\hat{y})+U_j(\hat{\mathbf{x}},\mu)=(\hat{\mathbf{x}},\hat{y})+
(U^1_j(\hat{\mathbf{x}},\mu),U^2_j(\hat{\mathbf{x}},\mu)),
\label{nfcv}
\end{equation}
(where $U^{1,2}_j$ are homogeneous degree $j$ polynomials in
the indicated variables, with coefficients respectively in
$\mathbb{R}^{2p+2}$ and $Q^1$)
system (\ref{y4}) can be put into formal normal form
\begin{subequations}\label{y5}
\begin{align}
\dot{\mathbf{x}}&={\displaystyle \mathring{B}\mathbf{x}+\sum_{j\geq
    2}\,\frac{1}{j!}\,g_j^1(\mathbf{x},y,\mu)}\label{y5:1}\\[0.15in]
{\displaystyle\frac{d}{dt}\,y}&={\displaystyle A_{Q^1}y+\sum_{j\geq
    2}\,\frac{1}{j!}\,g_j^2(\mathbf{x},y,\mu)}\label{y5:2}
\end{align}
\end{subequations}
such that the center manifold is locally given by $y=0$ and the local
flow of (\ref{y1}) on this center manifold is given by
\begin{equation}
\dot{\mathbf{x}}=\mathring{B}\mathbf{x}+\sum_{j\geq 2}\,\frac{1}{j!}\,g_j^1(\mathbf{x},0,\mu).
\label{y62}
\end{equation}
The nonlinear terms $g_j^1$ in (\ref{y62}) are in normal form in the
classical sense with respect to the matrix $\mathring{B}$.

\Section{Non-semisimple Equivariant Normal Form}

The matrix $\mathring{B}$ which appears in the previous section has a nilpotency
associated with it which somewhat complicates the computation of the
normal form (\ref{y62}).  In particular, the analysis which was
presented in \cite{CL1} in the semisimple case does not carry over here.  Nevertheless, in this
section we will show how to generalize the normal form analysis
presented in 
\cite{CL1}
to the present non-semisimple case.  Apart from solving the problem
which interests us in this paper, this approach may also shed some
light on the general case where Hypothesis \ref{spectralhyp} is
relaxed to include repeated eigenvalues with Jordan blocks.

\subsection{$\mathbb{T}^p\times\mathbb{R}$ normal forms}

Let ${\Psi}(0)$ denote the $(2p+1)\times 1$ matrix
obtained from the first $2p+1$ elements of the first column of
$\mathring{\Psi}(0)$ in Lemma \ref{psi0lem}, and let ${B}$ be the
$(2p+1)\times (2p+1)$ matrix
\begin{equation}
{B}=\mbox{\rm
  diag}(0,i\omega_1,-i\omega_1,\ldots,i\omega_p,-i\omega_p).
\label{B0def}
\end{equation}
It is easy to see from Lemma \ref{psi0lem} that $f_j^1$ in (\ref{fjdef}) is of the form
\[
f_j^1(\mathbf{x},0,\mu)=\left(\begin{array}{c}{\Psi}(0)\widehat{{\cal
        F}}_j(\mathring{\Phi}\,\mathbf{x},\mu)\\0\end{array}\right),
\]
and we will thus consider normal forms for the following
class of formal vector fields on $\mathbb{R}^{2p+2+s}$
\begin{equation}
\left(\begin{array}{c}\dot{x}\\\dot{\nu}\\\dot{\mu}\end{array}\right)=\left(\begin{array}{c}{B}x+\nu\,\mathbf{e}_{0}\\0\\0\end{array}\right)+
\sum_{j\geq
  2}\,\left(\begin{array}{c}f_j(x,\nu,\mu)\\0\\0\end{array}\right),
\label{y7ext}
\end{equation}
where we will use the convention that $\mathbf{e}_{k}$ is the unit row
vector in $\mathbb{R}^{q}$
with all entries equal to zero except the entry in the row $k+1$ whose
value is $1$, and where $q$ will depend on the context.

Using a mixture of complex and real coordinates,
we identify
\[
\mathbb{R}^{2p+2+s}=
\{(x_0,x_1,\overline{x_1},\ldots,x_p,\overline{x_p},\nu,\mu_1,\ldots,\mu_s)\,|\,x_0,\nu,\mu_j\in\mathbb{R},x_k\in\mathbb{C},\,j=1,\ldots,s,k=1,\ldots,p\}.
\]

Define the following $(2p+2+s)\times (2p+2+s)$ matrix
\[
\tilde{B}=\mbox{\rm diag}(\mathring{B},\mathbf{0}_{s}),
\] 
where $\mathbf{0}_s$ is the $s\times s$ zero matrix,
and let $\tilde{B}^T$ denote the
transpose of $\tilde{B}$.
Let 
\[
\Gamma=\overline{\{e^{s\tilde{B}^T}\,|\,s\in\mathbb{R}\}}
\]
where the
closure is taken in the space of $(2p+2+s)\times (2p+2+s)$ matrices.
Note
that $\Gamma$ is an abelian connected Lie group isomorphic
to
$\mathbb{T}^p\times\mathbb{R}$ where
$\mathbb{T}^p$ is the
$p$-torus:
\begin{equation}
\mathbb{T}^p=
\{\,\mbox{\rm
  diag}
(1,e^{i\theta_1},e^{-i\theta_1},\ldots,e^{i\theta_p},e^{-i\theta_p},
1,I_{s})\,\,\,|\,\,\,
\theta_j\in\mathbb{S}^1,\,j=1,\ldots,p\,\}
\label{Tpdef}
\end{equation}
and $\mathbb{R}$ is the one-parameter group parameterized as
\begin{equation}
\mathbb{R}=I_{2p+2+s}+\Theta\,\tilde{N}^T,\,\,\,\Theta\in\mathbb{R}
\end{equation}
where $\tilde{N}^T$ is the
transpose of the
nilpotent matrix 
\[
\tilde{N}=\tilde{B}-\mbox{\rm
  diag}(0,i\omega_1,-i\omega_1,\ldots,i\omega_p,-i\omega_p,0,\mathbf{0}_s).
\]
\begin{Def}
For a given integer $\ell\geq 2$ and a given normed space $X$, we denote
by $H^{2p+2+s}_{\ell}(X)$ the linear space of homogeneous polynomials
of degree $\ell$ in the $2p+2+s$ variables
$x=(x_0,x_1,\overline{x_1},\ldots,x_p,\overline{x_p})$, $\nu$ and
$\mu=(\mu_1,\ldots,\mu_s)$, with coefficients in $X$.
For $X=\mathbb{R}^{2p+2+s}$, define $H^{2p+2+s}_{\ell}(\mathbb{R}^{2p+2+s},\Gamma)\subset H^{2p+2+s}_{\ell}(\mathbb{R}^{2p+2+s})$
to be the
subspace of $\Gamma$-equivariant polynomials, i.e.
\[
\begin{array}{l}
\tilde{f}\in H^{2p+2+s}_{\ell}(\mathbb{R}^{2p+2+s},\Gamma)
\Longleftrightarrow\\[0.15in]
\tilde{f}\in H^{2p+2+s}_{\ell}(\mathbb{R}^{2p+2+s})\,\,\,\mbox{\rm and}\,\,\,
\gamma
\tilde{f}(\gamma^{-1}w)=\tilde{f}(w),\,\,\forall\,w=(x,\nu,\mu)\in\mathbb{R}^{2p+2+s},\,\,\forall\,\gamma\in\Gamma.
\end{array}
\]
\end{Def}

Normal forms for (\ref{y7ext}) are computed using the homological operator
\[
\begin{array}{c}
{\cal L}_{\tilde{B}}:H^{2p+2+s}_{\ell}(\mathbb{R}^{2p+2+s})\longrightarrow H^{2p+2+s}_{\ell}(\mathbb{R}^{2p+2+s})\\[0.15in]
\tilde{f}\longmapsto ({\cal L}_{\tilde{B}}\tilde{f})(w)=D\tilde{f}(w)\tilde{B}w-\tilde{B}\tilde{f}(w).
\end{array}
\]
To specify the normal form, we must find a complement in $H^{2p+2+s}_{\ell}(\mathbb{R}^{2p+2+s})$ to the
range of ${\cal L}_{\tilde{B}}$.  The following is a very well-known
result in the theory of normal forms \cite{ETBCI,GSSI,GSSII}.
\begin{prop}
\[
H^{2p+2+s}_{\ell}(\mathbb{R}^{2p+2+s})=H^{2p+2+s}_{\ell}(\mathbb{R}^{2p+2+s},\Gamma)\oplus\mbox{\rm
  range}\,{\cal L}_{\tilde{B}}
\]
\label{prop_enf1}
\end{prop}
It is
straightforward to compute the general element of $H^{2p+2+s}_{\ell}(\mathbb{R}^{2p+2+s},\Gamma)$.
\begin{lemma}
A smooth vector field
$\tilde{f}:\mathbb{R}^{2p+2+s}\longrightarrow\mathbb{R}^{2p+2+s}$
is $\mathbb{T}^p$-equivariant if and only if
$\tilde{f}$ is of the form
\begin{equation}
\tilde{f}(x,\nu,\mu)=
\left(
\begin{array}{c}
a_0(x_0,x_1\overline{x_1},\ldots,x_p\overline{x_p},\nu,\mu)\\ \\
a_1(x_0,x_1\overline{x_1},\ldots,x_p\overline{x_p},\nu,\mu)\,x_1\\ \\
\overline{a_1(x_0,x_1\overline{x_1},\ldots,x_p\overline{x_p},\nu,\mu)\,x_1}\\ \\
\vdots\\ \\
a_p(x_0,x_1\overline{x_1},\ldots,x_p\overline{x_p},\nu,\mu)\,x_p\\ \\
\overline{a_p(x_0,x_1\overline{x_1},\ldots,x_p\overline{x_p},\nu,\mu)\,x_p},\\
\\
b(x_0,x_1\overline{x_1},\ldots,x_p\overline{x_p},\nu,\mu)\\ \\
c_1(x_0,x_1\overline{x_1},\ldots,x_p\overline{x_p},\nu,\mu)\\ \\
\vdots\\ \\
c_s(x_0,x_1\overline{x_1},\ldots,x_p\overline{x_p},\nu,\mu),
\end{array}
\right)
\label{torus_nf}
\end{equation}
where $a_1,\ldots,a_p$ are smooth and complex-valued, and $a_{0}, b,
c_1,\ldots c_s$ are
smooth and real-valued.  
Furthermore, a vector field of the form (\ref{torus_nf}) is
$\Gamma$-equivariant (where $\Gamma\cong\mathbb{T}^p\times\mathbb{R}$)
if and only if
$a_0,a_1,\ldots,a_p,c_1,\ldots,c_{s}$ are
$\nu$-independent, and
\[
a_0=x_0\,g_0(x_0,x_1\overline{x_1},\ldots,x_p\overline{x_p},\mu),
\]
\[
b=\nu\,g_0(x_0,x_1\overline{x_1},\ldots,x_p\overline{x_p},\mu)+
g_1(x_0,x_1\overline{x_1},\ldots,x_p\overline{x_p},\mu)
\]
for some smooth functions $g_0$ and $g_1$.
\label{lemenf1}
\end{lemma}
\proof
The vector field\
$\tilde{f}$ is $\Gamma$-equivariant if and only if $\tilde{f}$ is both
$\mathbb{T}^p$-equivariant and $\mathbb{R}$-equivariant.
For the general form of the $\mathbb{T}^p$-equivariant vector field (\ref{torus_nf}),
see \cite{CL1}.  We then further require that (\ref{torus_nf})
commute with all matrices of the form 
$I_{2p+2+s}+\Theta\,\tilde{N}^T$, $\Theta\in\mathbb{R}$.  The result follows
after a straightforward computation.
\hfill\qed

The vector field $\tilde{f}$
in (\ref{y7ext}) has the special form $\tilde{f}=(f,0,0)$, which we want
our normal form changes of variables to preserve.  Since
we are only interested in the first $2p+1$ components of
(\ref{y7ext}), 
we would like to obtain a splitting of
$H^{2p+2+s}_{\ell}(\mathbb{R}^{2p+1})$ akin to the splitting 
of $H^{2p+2+s}_{\ell}(\mathbb{R}^{2p+2+s})$
in Proposition \ref{prop_enf1}.
For this purpose, we will need the following
\begin{Def}
\hspace*{1in}
\newline
\begin{enumerate}
\item[(a)]
We define $H^{2p+2+s}_{\ell}(\mathbb{R}^{2p+1},\mathbb{T}^p)$ to be the
subset of $H^{2p+2+s}_{\ell}(\mathbb{R}^{2p+1})$ consisting of
mappings
$f:\mathbb{R}^{2p+2+s}\longrightarrow\mathbb{R}^{2p+1}$
whose components are of the form of the first $2p+1$ components of
(\ref{torus_nf}).

As the notation suggests, note that
$H^{2p+2+s}_{\ell}(\mathbb{R}^{2p+1},\mathbb{T}^p)$ consists precisely
of elements of $H^{2p+2+s}_{\ell}(\mathbb{R}^{2p+1})$ which are
equivariant under an action of $\mathbb{T}^p$ on $\mathbb{R}^{2p+1}$:
\[
\begin{array}{l}
f\in
H^{2p+2+s}_{\ell}(\mathbb{R}^{2p+1},\mathbb{T}^p)\Longleftrightarrow\\ \\
f\in H^{2p+2+s}_{\ell}(\mathbb{R}^{2p+1})\,\,\,\mbox{\rm
  and}
\,\,\,f(\gamma_0\,x,\nu,\mu)=\gamma_0\,f(x,\nu,\mu),\,\,\,\forall\gamma_0\in\Gamma_0,\,\,\forall\,(x,\nu,\mu)\in\mathbb{R}^{2p+2+s},
\end{array}
\]
where $\Gamma_0$ is the group of $(2p+1)\times (2p+1)$ matrices which
is isomorphic to $\mathbb{T}^p$, and is parameterized as
\begin{equation}
\Gamma_0=
\{\,\mbox{\rm
  diag}(1,e^{i\theta_1},e^{-i\theta_1},\ldots,e^{i\theta_p},e^{-i\theta_p})\,\,\,|\,\,\,\theta_j\in\mathbb{S}^1,\,j=1,\ldots,p\,\}
\label{Tpdef2}
\end{equation}
\item[(b)] We define $H^{2p+1+s}_{\ell}(\mathbb{R}^{2p+1},\mathbb{T}^p)$ to
be the subspace of
$H^{2p+2+s}_{\ell}(\mathbb{R}^{2p+1},\mathbb{T}^p)$ consisting of
$\nu$-independent and $\mathbb{T}^p$-equivariant polynomials.
\item[(c)] 
We define the following operator
\[
\begin{array}{c}
{\cal L}_{{B},{\nu}} : H^{2p+2+s}_{\ell}(\mathbb{R}^{2p+1})
\longrightarrow H^{2p+2+s}_{\ell}(\mathbb{R}^{2p+1})\\[0.1in]
{\dps f\longmapsto ({\cal
  L}_{{B},\nu})(f)(x,\nu,\mu)=D_xf(x,\nu,\mu){B}x-{B}f(x,\nu,\mu)+
\nu\,\frac{\partial f}{\partial {x_0}}(x,\nu,\mu),}
\end{array}
\]
where ${B}$ is as in (\ref{B0def}),
and note that 
${\cal L}_{{B},\nu}$ is the usual homological operator
associated to the
$\dot{x}$ component of (\ref{y7ext}).  Furthermore, we have
${\cal L}_{\tilde{B}}(f,0,0)=({\cal L}_{{B},\nu}\,f,0,0)$.
\end{enumerate}
\end{Def}
\begin{prop}
\[
H^{2p+2+s}_{\ell}(\mathbb{R}^{2p+1})=H^{2p+1+s}_{\ell}(\mathbb{R}^{2p+1},\mathbb{T}^p)\oplus
\mbox{\rm range}\,{\cal L}_{{B},{\nu}}.
\]
\label{prop_enf2}
\end{prop}
\proof
The proof is given in the appendix.
\hfill
\qed

\subsection{Equivariant projection}
We will now construct an appropriate linear projection
associated with the splitting of
$H^{2p+2+s}_{\ell}(\mathbb{R}^{2p+1})$ given in Proposition
\ref{prop_enf2}.  
\begin{Def}
Let ${\displaystyle\int_{\Gamma_0}\,d\gamma}$ denote the normalized
Haar integral on $\Gamma_0\cong\mathbb{T}^p$ (see (\ref{Tpdef2})).  We define the linear
operator
\[
\begin{array}{c}
\mathring{A}:H^{2p+2+s}_{\ell}(\mathbb{R}^{2p+1})\longrightarrow H^{2p+2+s}_{\ell}(\mathbb{R}^{2p+1})\\[0.15in]
{\displaystyle f\longmapsto
(\mathring{A}\,f)(x,\nu,\mu)=\displaystyle\int_{\Gamma_0}\,\gamma\,f(\gamma^{-1}x,0,\mu)\,d\gamma}.
\end{array}
\]
\label{Adef}
\end{Def}
\begin{prop}
$\mathring{A}$ is a projection.  Furthermore,
\begin{equation}
\mbox{\rm range}\,\mathring{A}=
H^{2p+1+s}_{\ell}(\mathbb{R}^{2p+1},\mathbb{T}^p)
\label{Arange}
\end{equation}
and
\begin{equation}
\mbox{\rm ker}\,\mathring{A}=
\mbox{\rm range}\,{\cal L}_{{B},\nu}.
\label{Aker}
\end{equation}
\label{prop_Adecomposition}
\end{prop}
\proof
The proof is given in the appendix.
\hfill\qed

For any $f\in
H^{2p+2+s}_{\ell}(\mathbb{R}^{2p+1})$, write
\[
f=\mathring{A}f+(I-\mathring{A})f,
\]
and note that $\mathring{A}f$ is $\mathbb{T}^p$-equivariant and
$\nu$-independent and that $(I-\mathring{A})f\in
\mbox{\rm ker}\,\mathring{A}$.  
From Proposition \ref{prop_Adecomposition}, there exists
$h\in H^{2p+2+s}_{\ell}(\mathbb{R}^{2p+1})$ such that 
${\cal L}_{{B},\nu}h=(I-\mathring{A})f$.

\subsection{Phase decoupling}
Elements of the space
$H^{2p+1+s}_{\ell}(\mathbb{R}^{2p+1},\mathbb{T}^p)$ are
$\nu$-independent and equivariant
with respect to the torus group $\Gamma_0\cong\mathbb{T}^p$ defined in
(\ref{Tpdef2}).  As seen in \cite{CL1}, 
this toroidal equivariance can be used to
achieve a decoupling of the normal form.  Specifically, we have
\begin{prop}
Consider the following differential equation on $\mathbb{R}^{2p+1}$:
\[
\dot{x}={B}\,x+\nu\,\mathbf{e}_{0}+f(x,\mu),
\] 
where $f$ is smooth, $\nu$-independent,
satisfies $f(0,0)=0$, $Df(0,0)=0$, and is
$\Gamma_0\cong\mathbb{T}^p$-equivariant, i.e. $f$ has the form
\begin{equation}
f(x,\mu)=
\left(
\begin{array}{c}
a_0(x_0,x_1\overline{x_1},\ldots,x_p\overline{x_p},\mu)\\ \\
a_1(x_0,x_1\overline{x_1},\ldots,x_p\overline{x_p},\mu)\,x_1\\ \\
\overline{a_1(x_0,x_1\overline{x_1},\ldots,x_p\overline{x_p},\mu)\,x_1}\\ \\
\vdots\\ \\
a_p(x_0,x_1\overline{x_1},\ldots,x_p\overline{x_p},\mu)\,x_p\\ \\
\overline{a_p(x_0,x_1\overline{x_1},\ldots,x_p\overline{x_p},\mu)\,x_p},\end{array}
\right)
\label{torus_nf2}
\end{equation}
where $a_1,\ldots,a_p$ are smooth and complex-valued, and $a_{0}$ is
smooth and real-valued.  
Then under the under the change
of variables $x_0=\rho_0$, $x_j=\rho_je^{i\theta_j}$, $j=1,\ldots,p$, this
differential equation transforms into
\begin{equation}
\begin{array}{c}
\dot{\rho}_0=\nu+a_0(\rho_0,\rho_1^2,\ldots,\rho_p^2,\mu)\\ \\
\dot{\rho}_j=\mbox{\rm
  Re}(a_j(\rho_0,\rho_1^2,\ldots,\rho_p^2,\mu))\,\rho_j,\,\,\,j=1,\ldots,p\end{array}
\label{radial_eqs}
\end{equation}
and
\begin{equation}
\dot{\theta}_j=i\omega_j+
\mbox{\rm
  Im}(a_j(\rho_0,\rho_1^2,\ldots,\rho_p^2,\mu)),\,\,\,j=1,\ldots,p.
\label{angular_eqs}
\end{equation}
\end{prop}
\proof
This is a simple computation.
\hfill\qed

As in \cite{CL1},
we will call the subsystem (\ref{radial_eqs}) the {\em uncoupled
  radial part} of the normal form (\ref{torus_nf2}).  Recall that this uncoupled
radial part has some residual reflectional symmetry.
Denote by ${\mathbb Z}_{2,p}$ the group
whose action on $\mathbb{R}^{2p+1}$ is given by
\begin{equation}
(\rho_0,\rho_1,\ldots,\rho_p)\rightarrow(\rho_0,\lambda_1\rho_1,\ldots,\lambda_p\rho_p),
\label{multiz2act}
\end{equation}
where $\lambda_j\in\{1,-1\}$, $j=1,\ldots,p$.
\begin{Def}
For a given integer $\ell\geq 2$ and a given normed space $X$, we denote
by $H^{p+1+s}_{\ell}(X)$ the linear space of homogeneous polynomials
of degree $\ell$ in the $p+1+s$ variables 
$\rho=(\rho_0,\rho_1,\ldots,\rho_p)$ and
$\mu=(\mu_1,\ldots,\mu_s)$ with coefficients in $X$. 
Denote by
$H^{p+1+s}_{\ell}(\mathbb{R}^{p+1},\mathbb{Z}_{2,p})\subset H^{p+1+s}_{\ell}(\mathbb{R}^{p+1})$
the subspace of $H^{p+1+s}_{\ell}(\mathbb{R}^{p+1})$ 
consisting of
$\mathbb{Z}_{2,p}$-equivariant polynomials.
\end{Def}
It is easy to show (see \cite{GSSII}) that the most general element of
$H^{p+1+s}_{\ell}(\mathbb{R}^{p+1},\mathbb{Z}_{2,p})$ has the form of
the right-hand side of
(\ref{radial_eqs}).
We then define the following surjective linear mapping
\begin{equation}
\Pi : H^{2p+1+s}_{\ell}(\mathbb{R}^{2p+1},\mathbb{T}^p)\longrightarrow
H^{p+1+s}_{\ell}(\mathbb{R}^{p+1},\mathbb{Z}_{2,p})
\label{Pidef1}
\end{equation}
which is defined by sending the general element (\ref{torus_nf2}) of $H^{2p+1+s}_{\ell}(\mathbb{R}^{2p+1},\mathbb{T}^p)$
to the following element of $H^{p+1+s}_{\ell}(\mathbb{R}^{p+1},\mathbb{Z}_{2,p})$:
\begin{equation}
\left(
\begin{array}{c}
a_0(\rho_0,\rho_1^2,\ldots,\rho_p^2,\mu)\\
\mbox{\rm Re}(a_1(\rho_0,\rho_1^2,\ldots,\rho_p^2,\mu))\,\rho_1\\
\vdots\\
\mbox{\rm Re}(a_p(\rho_0,\rho_1^2,\ldots,\rho_p^2,\mu))\,\rho_p
\end{array}
\right).
\label{Pidef2}
\end{equation}

\subsection{Parameter splitting}
As seen in \cite{CL1}, it is useful when considering unfoldings
to be able to refine Propositions
\ref{prop_enf2} and \ref{prop_Adecomposition} in order to make
explicit the roles of $(x,\nu)$ as primary variables and $\mu$ as
unfolding parameters in (\ref{y5:1}).  For this purpose, we define the
following spaces.
\begin{Def}
Let $d\geq 1$ be a given integer (for our purposes, $d$ will be equal to
either $2p+1$ or to $p+1$), $\ell\geq 2$ be an integer, $X$ be a normed
linear space, and $G$ be a linear group acting on $X$.
Let $H^{d+1+s}_{\ell}(X)$ be the linear space of homogeneous
polynomials of degree $\ell$ in the variables $(\xi_1,\ldots,\xi_d)$,
$\nu$ and $\mu_1,\ldots,\mu_s$, and $H^{d+s}_{\ell}(X,G)\subset
H^{d+1+s}_{\ell}(X)$ be the subspace of $\nu$-independent and
$G$-equivariant polynomials.
\begin{enumerate}
\item[(a)]
We define $H^{d+1}_{\ell}(X)\subset H^{d+1+s}_{\ell}(X)$ and
$H^{d}_{\ell}(X,G)\subset H^{d+s}_{\ell}(X,G)$ to be the subspaces of
$\mu$-independent polynomials.
\item[(b)]
We define $P^{d+1+s}_{\ell}(X)\subset H^{d+1+s}_{\ell}(X)$ and
$P^{d+s}_{\ell}(X,G)\subset H^{d+s}_{\ell}(X,G)$ to be the subspaces
of polynomials which vanish at $\mu=0$.
\end{enumerate}
\end{Def}
It follows from these definitions that 
\begin{equation}
\begin{array}{c}
H^{2p+2+s}_{\ell}(\mathbb{R}^{2p+1})=H^{2p+2}_{\ell}(\mathbb{R}^{2p+1})\oplus
P^{2p+2+s}_{\ell}(\mathbb{R}^{2p+1})\\ \\
H^{2p+1+s}_{\ell}(\mathbb{R}^{2p+1},\mathbb{T}^p)=H^{2p+1}_{\ell}(\mathbb{R}^{2p+1},\mathbb{T}^p)\oplus
P^{2p+1+s}_{\ell}(\mathbb{R}^{2p+1},\mathbb{T}^p)\\ \\
H^{p+1+s}_{\ell}(\mathbb{R}^{p+1},\mathbb{Z}_{2,p})=H^{p+1}_{\ell}(\mathbb{R}^{p+1},\mathbb{Z}_{2,p})\oplus
P^{p+1+s}_{\ell}(\mathbb{R}^{p+1},\mathbb{Z}_{2,p}).
\end{array}
\label{paramnoparam}
\end{equation}
Furthermore, the various operators ${\cal L}_{{B},\nu}$, $\mathring{A}$
and $\Pi$ previously defined preserve these decompositions.  We then
get the following refinement of Propositions
\ref{prop_enf2} and \ref{prop_Adecomposition}:
\begin{prop}
\[
\begin{array}{ll}
{\cal L}_{{B},\nu}(H^{2p+2}_{\ell}(\mathbb{R}^{2p+1}))\subset
H^{2p+2}_{\ell}(\mathbb{R}^{2p+1}),\hspace*{0.3in} &
{\cal L}_{{B},\nu}(P^{2p+2+s}_{\ell}(\mathbb{R}^{2p+1}))\subset
P^{2p+2+s}_{\ell}(\mathbb{R}^{2p+1}),\\ \\
\mathring{A}(H^{2p+2}_{\ell}(\mathbb{R}^{2p+1}))=H^{2p+1}_{\ell}(\mathbb{R}^{2p+1},\mathbb{T}^p), &
\mathring{A}(P^{2p+2+s}_{\ell}(\mathbb{R}^{2p+1}))=P^{2p+1+s}_{\ell}(\mathbb{R}^{2p+1},\mathbb{T}^p),\\
\\
\Pi(H^{2p+1}_{\ell}(\mathbb{R}^{2p+1},\mathbb{T}^p))=H^{p+1}_{\ell}(\mathbb{R}^{p+1},\mathbb{Z}_{2,p}), &
\Pi(P^{2p+1+s}_{\ell}(\mathbb{R}^{2p+1},\mathbb{T}^p))=P^{p+1+s}_{\ell}(\mathbb{R}^{p+1},\mathbb{Z}_{2,p})
\end{array}
\]
Consequently, if we define $\mathring{A}|_1$ and $\mathring{A}|_2$ to be respectively the
restrictions of $\mathring{A}$ on $H^{2p+2}_{\ell}(\mathbb{R}^{2p+1})$ and on
$P^{2p+2+s}_{\ell}(\mathbb{R}^{2p+1})$, 
then $\mathring{A}|_1$ and $\mathring{A}|_2$ are projections.
Similarly define
${\cal L}_{{B},\nu}|_1$ and ${\cal L}_{{B},\nu}|_2$
to be respectively the
restrictions of ${\cal L}_{{B},\nu}$ on $H^{2p+2}_{\ell}(\mathbb{R}^{2p+1})$ and on
$P^{2p+2+s}_{\ell}(\mathbb{R}^{2p+1})$, then
\[
\begin{array}{ll}
\mbox{\rm
  range}\,\mathring{A}|_1=H^{2p+1}_{\ell}(\mathbb{R}^{2p+1},\mathbb{T}^p),\hspace*{0.3in}&
\mbox{\rm
  range}\,\mathring{A}|_2=P^{2p+1+s}_{\ell}(\mathbb{R}^{2p+1},\mathbb{T}^p),\\ \\
\mbox{\rm ker}\,\mathring{A}|_1=\mbox{\rm range}\,{\cal
  L}_{{B},\nu}|_1,&
\mbox{\rm ker}\,\mathring{A}|_2=\mbox{\rm range}\,{\cal
  L}_{{B},\nu}|_2,
\end{array}
\]
\[
\begin{array}{c}
H^{2p+2}_{\ell}(\mathbb{R}^{2p+1})=H^{2p+1}_{\ell}(\mathbb{R}^{2p+1},\mathbb{T}^p)\oplus\mbox{\rm
  range}\,{\cal L}_{{B},\nu}|_1,\\ \\
P^{2p+2+s}_{\ell}(\mathbb{R}^{2p+1})=P^{2p+1+s}_{\ell}(\mathbb{R}^{2p+1},\mathbb{T}^p)\oplus\mbox{\rm
  range}\,{\cal L}_{{B},\nu}|_2.
\end{array}
\]
\label{prop_enf3}
\end{prop}

We now combine the results of this section with the Faria and
Magalh$\tilde{\mbox{\rm a}}$es normal form procedure outlined in
Section 2 in order to obtain the following version of Theorem 5.8 of
\cite{FMTB} and Theorem 2.16 of \cite{FMH}.
\begin{thm}
Consider the system (\ref{y4}) 
\begin{equation}
\begin{array}{rcl}
\dot{\mathbf{x}}&=&{\dps \mathring{B}\,\mathbf{x}+\sum_{j\geq
    2}\,f_j^1(\mathbf{x},y,\mu)}\\[0.15in]
{\dps\frac{d}{dt}y}&=&{\dps A_{Q^1}\,y+\sum_{j\geq
    2}\,f_j^2(\mathbf{x},y,\mu),}
\end{array}
\label{y7tayl}
\end{equation}
where ${\dps\mathbf{x}=\left(\begin{array}{c}x\\\nu\end{array}\right)}$.
Write
\begin{equation}
f_j^1(\mathbf{x},0,\mu)=\left(\begin{array}{c}{\Psi}(0)\widehat{{\cal
        F}}_j(\mathring{\Phi}\,\mathbf{x},\mu)\\0\end{array}\right)=\left(\begin{array}{c}
h_j(\mathbf{x})+q_j(\mathbf{x},\mu)\\ 0\end{array}\right),
\label{HPdecomfj}
\end{equation}
where $h_j\in H^{2p+2}_{\ell}(\mathbb{R}^{2p+1})$ and $q_j\in
P^{2p+2+s}_{\ell}(\mathbb{R}^{2p+1})$.
Then there is a formal near-identity change of variables
\[
(\mathbf{x},y)\rightarrow
(\hat{\mathbf{x}},\hat{y})+(U^1(\mathbf{x}),U^2(\mathbf{x}))+(W^1(\mathbf{x},\mu),W^2(\mathbf{x},\mu))
\]
(where $W^1(\mathbf{x},0)=0$ and $W^2(\mathbf{x},0)=0$) which
transforms (\ref{y7tayl}) into system (\ref{y5}) (upon dropping the
hats), and the flow on the invariant local center manifold $y=0$ is
given by
\begin{equation}\label{y6}
\begin{array}{rcl}
\dot{x}&=&{\dps {B}\,x+\nu\,\mathbf{e}_{0}+\sum_{j\geq
  2}\,((\mathring{A}|_1(h_j+Y_j))(x)+(\mathring{A}|_2(q_j+Z_j))(x,\mu)),}\\[0.15in]
\dot{\nu}&=&0\\[0.15in]
\dot{\mu}&=&0,
\end{array}
\end{equation}
where $Y_2=0$, $Z_2=0$, and for $j\geq 3$, $Y_j=Y_j(x,\nu)$ and $Z_j=
Z_j(x,\nu,\mu)$
are the extra contributions to the terms of order $j$ coming from the
transformation of the lower order ($<j$) terms, and $Z_j$ vanishes at $\mu=0$.
\label{thmdsdecomnf}
\end{thm}

\Section{Main results}
In the semisimple case of \cite{CL1}, the main realizability results are a
consequence of the surjectivity of a certain linear operator between
suitable spaces of polynomials, which is proven in Proposition 4.3 of
\cite{CL1}.
In this section, we define the corresponding linear operator in the
present non-semisimple case, and prove that its surjectivity follows
from the surjectivity in the semisimple case.  Once this
has been achieved, we will state our main realizability results for the
class of uncoupled radial equations (\ref{Pidef2}) which are obtained
from writing the normal form center manifold equations (\ref{y6}) in
polar coordinates.

\subsection{Linear analysis}

\begin{Def}
For a given integer $\ell\geq 2$,
let $H^{2p+2+s}_{\ell}(\mathbb{R})$ denote the linear space of
homogeneous degree $\ell$ polynomials in the $2p+2$ variables
\[
\mathbf{v}=(\mathbf{v}_0,\mathbf{v}_1,\ldots,\mathbf{v}_p)=\left(\,\left(\begin{array}{c}v_0\\w_0\end{array}\right)\,,\,
\left(\begin{array}{c}v_1\\w_1\end{array}\right)\,,\,\ldots\,,\,
\left(\begin{array}{c}v_p\\w_p\end{array}\right)\,\right),
\]
and $s$ parameters $\mu=(\mu_1,\ldots,\mu_s)$ with real coefficients.
Denote by $H^{p+1+s}_{\ell}(\mathbb{R})$ the linear space of
homogeneous degree $\ell$ polynomials in the variables
$v=(v_0,v_1,\ldots,v_p)$ and $\mu=(\mu_1,\ldots,\mu_s)$.
Note that we may decompose these spaces as in (\ref{paramnoparam}) as follows
\[
H^{2p+2+s}_{\ell}(\mathbb{R})=H^{2p+2}_{\ell}(\mathbb{R})\oplus
P^{2p+2+s}_{\ell}(\mathbb{R})
\]
and
\[
H^{p+1+s}_{\ell}(\mathbb{R})=H^{p+1}_{\ell}(\mathbb{R})\oplus
P^{p+1+s}_{\ell}(\mathbb{R})
\]
where $H^{2p+2}_{\ell}(\mathbb{R})$ and $H^{p+1}_{\ell}(\mathbb{R})$
are the $\mu$-independent polynomials, and
$P^{2p+2+s}_{\ell}(\mathbb{R})$ and $P^{p+1+s}_{\ell}(\mathbb{R})$ are
the polynomials which vanish at $\mu=0$.
Finally, define the surjective linear mapping
\[
{\cal R}:H^{2p+2+s}_{\ell}(\mathbb{R})\longrightarrow
H^{p+1+s}_{\ell}(\mathbb{R})
\]
as
\begin{equation}
({\cal R}(h))(v,\mu)=h\left(\,\left(\begin{array}{c}v_0\\0\end{array}\right)\,,\,\left(\begin{array}{c}v_1\\0\end{array}\right)\,,\,\ldots\,,\,
\left(\begin{array}{c}v_p\\0\end{array}\right)\,,\,\mu\right).
\label{Rprojdef}
\end{equation}
\label{Vspacedef}
\end{Def}

Let $\mathring{\Phi}$ be as in (\ref{phichoice}).  
We will define ${\Phi}$ to be the $1\times (2p+1)$ matrix
obtained from the first $2p+1$ elements of the first row of $\mathring{\Phi}$,
i.e.
\[
{\Phi}=(\,1\,\,\,e^{i\omega_1\theta}\,\,\,e^{-i\omega_1\theta}\,\,\,\cdots\,\,\,e^{i\omega_p\theta}\,\,\,e^{-i\omega_p\theta}\,).
\]
From
Lemma \ref{psi0lem} and the fact that
 ${\Psi}(0)$ denotes the $(2p+1)\times 1$ matrix
obtained from the first $2p+1$ elements of the first column of
$\mathring{\Psi}(0)$ in Lemma \ref{psi0lem}, it follows that
\begin{equation}
{\Psi}(0)=
\mbox{\rm
    col}(u_0,u_1,\overline{u_1},\ldots,u_p,\overline{u_p})
\label{Vdef}
\end{equation}
where $u_0\neq 0$ is real and $u_j\neq 0$ are complex, $j=1,\ldots,p$.

Let $\tau=(\tau_1,\ldots,\tau_{p+1})\in\mathbb{R}^{p+1}$, and define
\begin{equation}
\mathring{E}_{\tau}=\left(\begin{array}{c}\mathring{\Phi}(\tau_1)\\\vdots\\\mathring{\Phi}(\tau_d)\end{array}\right)=
\left(\begin{array}{ccccccc}
1&e^{i\omega_1\tau_1}&e^{-i\omega_1\tau_1}&\cdots&e^{i\omega_p\tau_1}&e^{-i\omega_p\tau_1}&\tau_1\\
0&0&0&\cdots&0&0&1\\
1&e^{i\omega_1\tau_2}&e^{-i\omega_1\tau_2}&\cdots&e^{i\omega_p\tau_2}&e^{-i\omega_p\tau_2}&\tau_2\\
0&0&0&\cdots&0&0&1\\
\vdots&\vdots&\vdots&\vdots&\vdots&\vdots&\vdots\\
1&e^{i\omega_1\tau_{p+1}}&e^{-i\omega_1\tau_{p+1}}&\cdots&e^{i\omega_p\tau_{p+1}}&e^{-i\omega_p\tau_{p+1}}&\tau_{p+1}\\
0&0&0&\cdots&0&0&1
\end{array}\right)
\label{Edefmulthopf}
\end{equation}
\begin{equation}
{E}_{\tau}=\left(\begin{array}{c}{\Phi}(\tau_1)\\\vdots\\{\Phi}(\tau_d)\end{array}\right)=
\left(\begin{array}{cccccc}
1&e^{i\omega_1\tau_1}&e^{-i\omega_1\tau_1}&\cdots&e^{i\omega_p\tau_1}&e^{-i\omega_p\tau_1}\\
1&e^{i\omega_1\tau_2}&e^{-i\omega_1\tau_2}&\cdots&e^{i\omega_p\tau_2}&e^{-i\omega_p\tau_2}\\
\vdots&\vdots&\vdots&\vdots&\vdots&\vdots\\
1&e^{i\omega_1\tau_{p+1}}&e^{-i\omega_1\tau_{p+1}}&\cdots&e^{i\omega_p\tau_{p+1}}&e^{-i\omega_p\tau_{p+1}}
\end{array}\right).
\label{Edefmulthopf0}
\end{equation}
As in \cite{CL1}, we define the $\ell$-mappings associated to $\mathring{E}_{\tau}$
and to ${E}_{\tau}$ respectively as follows:
\begin{equation}
\begin{array}{c}
\mathring{\cal E}^{\ell}_{{\tau}}:H^{2p+2+s}_{\ell}(\mathbb{R})\longrightarrow
H^{2p+2+s}_{\ell}(\mathbb{R}^{2p+1})\\[0.1in]
(\mathring{\cal
  E}^{\ell}_{\tau}(h))(\mathbf{x},\mu)\equiv
{\Psi}(0)\,h(\mathring{E}_{\tau}\,\mathbf{x},\mu)
\end{array}
\label{Ecaldefmulthopf}
\end{equation}
and
\begin{equation}
\begin{array}{c}
{\cal E}^{\ell}_{{\tau}}:H^{p+1+s}_{\ell}(\mathbb{R})\longrightarrow
H^{2p+1+s}_{\ell}(\mathbb{R}^{2p+1})\\[0.1in]
{\cal
  E}^{\ell}_{\tau}(h)(x,\mu)\equiv
{\Psi}(0)\,h({E}_{\tau}\,{x},\mu)
\end{array}
\label{Ecaldefmulthopf0}
\end{equation}

Now, let
$\Pi:H^{2p+1+s}_{\ell}(\mathbb{R}^{2p+1},\mathbb{T}^p)\longrightarrow H^{p+1+s}_{\ell}(\mathbb{R}^{p+1},\mathbb{Z}_{2,p})$
be the mapping defined in (\ref{Pidef1})-(\ref{Pidef2}), let
$\mathring{A}:H^{2p+2+s}_{\ell}(\mathbb{R}^{2p+1})\longrightarrow
H^{2p+1+s}_{\ell}(\mathbb{R}^{2p+1},\mathbb{T}^p)$ be the projection
operator defined in
Definition \ref{Adef}, and define
\[
{A}:H^{2p+1+s}_{\ell}(\mathbb{R}^{2p+1})\longrightarrow
H^{2p+1+s}_{\ell}(\mathbb{R}^{2p+1},\mathbb{T}^p)
\]
by
\[
{A}g(x,\mu)=\int_{\Gamma_0}\,\gamma\,g(\gamma^{-1}\,x,\mu)\,d\gamma,
\]
as in \cite{CL1}.
Our main result in this section is the following:
\begin{prop}
For an open and dense set ${\cal U}\subset\mathbb{R}^{p+1}$,
the following linear mapping is surjective for all $\tau\in {\cal U}$:
\[
\Pi\circ \mathring{A}\circ \mathring{\cal E}^{\ell}_{{\tau}} : H^{2p+2+s}_{\ell}(\mathbb{R})\longrightarrow
H^{p+1+s}_{\ell}(\mathbb{R}^{p+1},\mathbb{Z}_{2,p}).
\]
\label{MultipleHopfinv}
\end{prop}
\proof
It is a simple computation to show that
\[
\Pi\circ \mathring{A}\circ \mathring{\cal E}^{\ell}_{\tau}=
(\Pi\circ {A}\circ {\cal E}^{\ell}_{\tau})\circ {\cal R},
\]
where ${\cal R}$ is the surjective linear operator defined in (\ref{Rprojdef}).
In \cite{CL1}, it was shown that there exists an open and dense set ${\cal
  U}\subset\mathbb{R}^{p+1}$ such that for all $\tau\in {\cal U}$, the mapping
\[
\Pi\circ {A}\circ {\cal E}^{\ell}_{\tau} :
H^{p+1+s}_{\ell}(\mathbb{R})\longrightarrow
H^{p+1+s}_{\ell}(\mathbb{R}^{p+1},\mathbb{Z}_{2,p})
\]
is surjective.
Thus, for all $\tau\in {\cal U}$, $\Pi\circ \mathring{A}\circ \mathring{\cal
  E}^{\ell}_{\tau}$ is surjective.
\hfill
\qed

\subsection{Main results on realizability}
With Proposition \ref{MultipleHopfinv} in hand, we obtain
realizability results analogous to the semisimple case.  Since the
proofs are almost identical, we merely give the statements of these
results here.  First, we will
define the following linear spaces of (non-homogeneous) polynomials
\begin{Def}
\[
\begin{array}{ll}
{\cal H
  }^{p+1+s}_{\ell}(\mathbb{R})\equiv\oplus_{j=2}^{\ell}\,H^{p+1+s}_{j}(\mathbb{R}),\hspace*{1in}&
{\cal
  H}^{p+1}_{\ell}(\mathbb{R})\equiv\oplus_{j=2}^{\ell}\,H^{p+1}_j(\mathbb{R}),\\[0.15in]
{\cal P}^{p+1+s}_{\ell}(\mathbb{R})\equiv\oplus_{j=2}^{\ell}\,P^{p+1+s}_j(\mathbb{R}),&
{\cal
  H}^{p+1+s}_{\ell}(\mathbb{R}^{p+1},\mathbb{Z}_{2,p})\equiv\oplus_{j=2}^{\ell}\,H^{p+1+s}_{j}(\mathbb{R}^{p+1},\mathbb{Z}_{2,p}),\\[0.15in]
{\cal
  H}^{p+1}_{\ell}(\mathbb{R}^{p+1},\mathbb{Z}_{2,p})\equiv\oplus_{j=2}^{\ell}\,H^{p+1}_j(\mathbb{R}^{p+1},\mathbb{Z}_{2,p}),&
{\cal P}^{p+1+s}_{\ell}(\mathbb{R}^{p+1},\mathbb{Z}_{2,p})\equiv\oplus_{j=2}^{\ell}\,P^{p+1+s}_{j}(\mathbb{R}^{p+1},\mathbb{Z}_{2,p}).
\end{array}
\]
\label{nonhompolys}
\end{Def}
\begin{thm}
Consider the RFDE (\ref{y1split}), and let $\Lambda_0$ denote the set of
solutions of (\ref{ychar}) with zero real part.  Suppose that
Hypothesis \ref{spectralhyp} is satisfied. 
Let $\ell\geq 2$ be a given integer.  For each
$h\in {\cal H}^{p+1}_{\ell}(\mathbb{R}^{p+1},\mathbb{Z}_{2,p})$
and each
$q\in {\cal P}^{p+1+s}_{\ell}(\mathbb{R}^{p+1},\mathbb{Z}_{2,p})$
there are $p+1$ distinct points
$\tau_1,\ldots,\tau_{p+1}\in [-r,0]$, an $\eta\in {\cal
  H}^{p+1}_{\ell}(\mathbb{R})$
and a $\xi\in {\cal P}^{p+1+s}_{\ell}(\mathbb{R})$
such that if 
\[
\hat{F}(z_t,\nu,\mu)=\eta(z(t+\tau_1),\ldots,z(t+\tau_d))+\xi(z(t+\tau_1),\ldots,z(t+\tau_d),\mu)
\]
in
  (\ref{y1split}), then in polar coordinates, the radial part of the center manifold equations
  (\ref{y6}) in $\mathbb{T}^p\times\mathbb{R}$-equivariant normal form up to degree $\ell$ reduces to
  $\dot{\rho}=\nu\,\mathbf{e}_0+h(\rho)+q(\rho,\mu)$, where
  $\rho\equiv (\rho_0,\rho_1,\ldots,\rho_p)$.
  In fact, $\tau\equiv (\tau_1,\ldots,\tau_{p+1})$ can be chosen in an open and
  dense set of $[-r,0]^{p+1}$, independently of the particular $h$ and $q$
  to be
  realized (i.e. only $\eta$ and $\xi$ must be changed in order to account for
  different jets to be realized).
\label{mainthm1}
\end{thm}

As in the semisimple case of \cite{CL1}, we can show that the number of
delays $p+1$ shown above to be sufficient to solve the realizability
problem for the radial part is optimal i.e. surjectivity is violated
beyond some finite order $\ell_0$
if the number of delays is less than $p+1$ (see Theorem 5.4 of \cite{CL1}).

For general RFDEs,
we have the following result on realization of unfoldings:
\begin{thm}
Consider the general nonlinear RFDE
\begin{equation}
\dot{z}(t)=L_0\,z_t+\nu+N(z_t)
\label{genunpertrfde}
\end{equation}
where $\nu\in\mathbb{R}$, $L_0:C_1\rightarrow\mathbb{R}$ is
a bounded linear
operator from $C_1\equiv C\left(  \left[  -r,0\right]  ,\mathbb{R}\right)$ into $\mathbb{R}$, and $N$ is
a smooth function from $C_1$ into $\mathbb{R}$,
with $N(0)=0$, 
$DN(0)=0$.
Let $\Lambda_0$ denote the set of
solutions of (\ref{ychar}) with zero real part and suppose that
Hypothesis \ref{spectralhyp} is satisfied. 
Then the local dynamics of (\ref{genunpertrfde})
near the origin on an invariant
center manifold can be described by a system of ordinary differential
equations on $\mathbb{R}^{2p+1}$.  Moreover, this ODE system can be
brought into $\mathbb{T}^p\times\mathbb{R}$-equivariant normal form to any desired
order $\ell$, and the resulting (truncated at order $\ell$) normal
form can be uncoupled into an $p+1$-dimensional system and a
$p$-dimensional system
\begin{eqnarray}
\dot{\rho}&=&\nu\,\mathbf{e}_0+h(\rho\,;\,N)\label{1st}\\[0.15in]
\dot{\theta}&=&k(\rho\,;\,N),\label{2nd}
\end{eqnarray}
where for given $N$, $h(\cdot\,;\,N)$ is some element of ${\cal
  H}^{p+1}_{\ell}(\mathbb{R}^{p+1},\mathbb{Z}_{2,p})$, and
  $k(\cdot\,;\,N):\mathbb{R}^p\longrightarrow\mathbb{R}^p$.
Let $\tilde{h}(\rho,\mu)$ be an $s$-parameter equivariant unfolding of $h$ of degree at
  most $\ell$, i.e.
$\tilde{h}\in {\cal H}^{p+1+s}_{\ell}(\mathbb{R}^{p+1},\mathbb{Z}_{2,p})$ and
  $\tilde{h}(\cdot,0)=h(\cdot\,;\,N)$.
Then there exists an $s$-parameter unfolding of (\ref{genunpertrfde})
  of the form
\begin{equation}
\dot{z}(t)=L_0(z_t)+\nu+N(z_t)+\xi(z(t+\tau_1),\ldots,z(t+\tau_{p+1}),\mu)
\label{genunfrfde}
\end{equation}
(where $\tau=(\tau_1,\ldots,\tau_{p+1})\in\mathbb{R}^{p+1}$, and
$\xi\in {\cal P}^{p+1+s}_{\ell}(\mathbb{R})$ vanishes at $\mu=0$)
which realizes the unfolded radial equations
\[
\dot{\rho}=\nu\,\mathbf{e}_0+\tilde{h}(\rho,\mu)
\]
on an invariant center manifold for (\ref{genunfrfde}).
\label{mainthm2}
\end{thm}

\subsection{Generic saddle-node/Hopf interaction}

The following example is an immediate consequence of our results.

\begin{examp}
{\rm
Consider the RFDE (\ref{y1split}) in the case $\mu=0$,
\begin{equation}
\dot{z}(t)=L_0z_t+\nu+\hat{F}(z_t),
\label{y50}
\end{equation}
such that the
characteristic equation (\ref{ychar}) has simple purely imaginary
roots $\pm i\omega\neq 0$, a simple root at $0$, and no other
roots on the imaginary axis.
If
\begin{equation}
\begin{array}{lll}
\hat{F}(z_t)&=&A_{20}(z(t+\tau_1))^2+A_{11}z(t+\tau_1)z(t+\tau_2)+A_{02}(z(t+\tau_2))^2\\
&&A_{30}(z(t+\tau_1))^3+A_{21}z((t+\tau_1))^2z(t+\tau_2)+A_{12}z(t+\tau_1)(z(t+\tau_2))^2+A_{03}(z(t+\tau_2))^3,
\end{array}
\label{FsshopfFM}
\end{equation}
where $\tau_1,\tau_2\in [-r,0]$,
then the
uncoupled radial part of the center manifold equations to cubic order
are the following Guckenheimer \cite{Guck1,Guck2} normal form
\begin{equation}
\begin{array}{lll}
\dot{\rho_0}&=&\nu+a_{1}\rho_0^2+a_{2}\rho_1^2+a_3\rho_0^3+a_4\rho_0\rho_1^2\\[0.15in]
\dot{\rho_1}&=&b_1\rho_0\rho_1+b_2\rho_1^3+b_3\rho_1\rho_0^2,
\label{sshopfnfFM}
\end{array}
\end{equation}
where the coefficients $a_{1,2,3,4}$ and $b_{1,2,3}$ are functions of
$A_{20}$, $A_{11}$, $A_{02}$, $A_{30}$, $A_{21}$, $A_{12}$, $A_{03}$,
$\tau_1$ and $\tau_2$.  From our results, it follows that generically, any values of $a_{1,2,3,4}$ and
$b_{1,2,3}$ can be achieved by appropriate choice of
$A_{20}$, $A_{11}$, $A_{02}$, $A_{30}$, $A_{21}$, $A_{12}$, $A_{03}$.
Also, the following versal unfolding
of (\ref{sshopfnfFM})
\begin{equation}
\begin{array}{lll}
\dot{\rho_0}&=&\nu+a_{1}\rho_0^2+a_{2}\rho_1^2+a_3\rho_0^3+a_4\rho_0\rho_1^2\\[0.15in]
\dot{\rho_1}&=&\mu\rho_1+b_1\rho_0\rho_1+b_2\rho_1^3+b_3\rho_1\rho_0^2,
\label{sshopfnfFMunf}
\end{array}
\end{equation}
is generically realized by the
following unfolding of (\ref{y50})
\[
\dot{z}(t)=L_0z_t+\nu+\hat{F}(z_t)+\mu\,(A_{10}z(t+\tau_1)+A_{01}\mu_2\,z(t+\tau_2))
\]
for appropriate choice of $A_{10}$ and $A_{01}$.}
\end{examp}

\Section{Concluding remarks}
Our solution here to the realizability problem for the non-semisimple
saddle-node/multiple Hopf interaction in scalar RFDEs
complements our results in \cite{CL1} for the semisimple cases of
transcritical/multiple Hopf interaction and non-resonant multiple Hopf
bifurcation.  We note that our results here on realizability are
generic,
optimal in the number of delays required to guarantee realizability,
and applicable to any finite order expansion and truncation of the
normal
form. Therefore, nonlinear degeneracies
and their unfoldings for the saddle-node/multiple Hopf interaction are
covered by our theory.

As in \cite{CL1}, we have not considered general $n>1$ dimensional systems
of RFDEs.  This is not based on any deep theoretical issues associated
to the $n>1$ case,
but rather to complications arising out of notation and messy
algebraic computations which would make the exposition extremely
cumbersome.
We have, however, every reason to believe that the realizability
problem for each of the bifurcations studied in \cite{CL1} and in the
present paper could be solved for systems using the techniques and framework
we have developed. 

More subtle is the issue of relaxing Hypothesis \ref{spectralhyp} to
include repeated eigenvalues with Jordan blocks, and rational
resonances in the purely imaginary eigenvalues.  Apart from affecting
the dimension of the torus group admitted by the normal form (and
consequently
the dimension of the uncoupled radial equations),
Jordan blocks would introduce an additional non-compact component to
the normal form symmetry.  In this case, our present analysis
in this paper may shed some valuable light on a suitable approach to tackle
the general problem, i.e. first solve the semisimple case, and then
try to exploit the semisimple solution as much as possible in order
to prove realizability in the associated non-semisimple problem.

\vspace*{0.25in}
\noindent
{\Large\bf Acknowledgments}

\vspace*{0.2in}
This research is partly supported by the
Natural Sciences and Engineering Research Council of Canada in the
form of a Discovery Grant (VGL), by a Premier's Research Excellence
Award from the Ontario Ministry of Economic Development and Trade and
the University of Ottawa (VGL), and by the
Centre de Recherches Math\'ematiques (YC).

\appendix

\Section{Proof of Proposition~\ref{prop_enf2}}

Let $f$ be a given element of
$H^{2p+2+s}_{\ell}(\mathbb{R}^{2p+1})$, and consider
$\tilde{f}=(f,0,0)\in H^{2p+2+s}_{\ell}(\mathbb{R}^{2p+2+s})$.  From
Proposition \ref{prop_enf1}, there exists $h=(h^x,h^{\nu},h^{\mu})\in
H^{2p+2+s}_{\ell}(\mathbb{R}^{2p+2+s})$ and a unique
$g=(g^x,g^{\nu},g^{\mu})\in
H^{2p+2+s}_{\ell}(\mathbb{R}^{2p+2+s},\Gamma)$ such that
\begin{equation}
(f,0,0)={\cal L}_{\tilde{B}}(h^x,h^{\nu},h^{\mu})+(g^x,g^{\nu},g^{\mu}).
\label{hom_extend}
\end{equation}
Note that
\[
{\cal
  L}_{\tilde{B}}(h^x,h^{\nu},h^{\mu})=\left({\cal
  L}_{{B},\nu}\,h^x-h^{\nu}\,\mathbf{e}_{0}\,,\,D_xh^{\nu}{B}x+\nu\frac{\partial
  h^{\nu}}{\partial x_0}\,,\,D_xh^{\mu}{B}x+\nu\frac{\partial
  h^{\mu}}{\partial x_0}\right).
\]
Thus, from (\ref{hom_extend}), we get
\begin{equation}
D_xh^{\nu,\mu}{B}x+\nu\frac{\partial h^{\nu,\mu}}{\partial {x_0}}+g^{\nu,\mu}=0
\label{2ndcompnt}
\end{equation}
and
\begin{equation}
f={\cal L}_{{B},\nu}h^{x}-h^{\nu}\,\mbox{\bf e}_{0}+g^x.
\label{1stcompnt}
\end{equation}
From Lemma \ref{lemenf1}, 
$g^{\nu}$ is of the form
\[
g^{\nu}=\nu\,r_{1}(x_0,x_1\overline{x_1},\ldots,x_p\overline{x_p},\mu)+
r_2(x_0,x_1\overline{x_1},\ldots,x_p\overline{x_p},\mu),
\]
where $r_1$ is such that the component of
$g^x$ along $\mathbf{e}_0$, $g^x_0$, is of the form
\[
g^x_{0}=x_0\,r_1(x_0,x_1\overline{x_1},\ldots,x_p\overline{x_p},\mu).
\]
Therefore, the $\nu$ part of equation
(\ref{2ndcompnt}) reduces to
\begin{equation}
\nu\,\frac{\partial h^{\nu}}{\partial x_0}+\sum_{j=1}^p\,i\omega_j
\left(x_j\frac{\partial h^{\nu}}{\partial
    {x_j}}-\overline{x_j}\frac{\partial h^{\nu}}{\partial {\overline{x_j}}}\right)+\nu\,r_1+r_2=0.
\label{2ndcompnts}
\end{equation}
Before we proceed any further, we will need the following two lemmas.
\begin{lemma}
Let $h$ be a real-valued smooth function of
$x_0,x_1,\overline{x_1},\ldots,x_p,\overline{x_p},\nu$ and $\mu$ such
that the function
$\mathfrak{g}(x_0,x_1,\overline{x_1},\ldots,x_p,\overline{x_p},\nu,\mu)$ defined
by
\begin{equation}
\mathfrak{g}(x_0,x_1,\overline{x_1},\ldots,x_p,\overline{x_p},\nu,\mu)\equiv
\sum_{j=1}^p\,i\omega_j\left(x_j\frac{\partial h}{\partial
    x_j}-\overline{x_j}\frac{\partial
    h}{\partial\overline{x_j}}\right)
\label{wassuppp1}
\end{equation}
is $\mathbb{T}^p$ invariant, i.e.
\begin{equation}
\mathfrak{g}(x_0,e^{i\theta_1}x_1,e^{-i\theta_1}\overline{x_1},\ldots,e^{i\theta_p}x_p,e^{-i\theta_p}\overline{x_p},\nu,\mu)=
\mathfrak{g}(x_0,x_1,\overline{x_1},\ldots,x_p,\overline{x_p},\nu,\mu)
\label{wassuppp2}
\end{equation}
for all $\theta_1,\ldots,\theta_p\in\mathbb{R}$, and for all
$(x,\nu,\mu)\in\mathbb{R}^{2p+2+s}$.
Then $h$ is also $\mathbb{T}^p$ invariant, and $\mathfrak{g}=0$.
\label{applem1}
\end{lemma}
\proof
A simple computation using (\ref{wassuppp1}) and (\ref{wassuppp2})
leads to
\[
\frac{d}{d\xi}
h(x_0,e^{i\omega_1\xi}x_1,e^{-i\omega_1\xi}\overline{x_1},\ldots,e^{i\omega_p\xi}x_p,e^{-i\omega_p\xi}\overline{x_p},\nu,\mu)=
\mathfrak{g}(x_0,x_1,\overline{x_1},\ldots,x_p,\overline{x_p},\nu,\mu).
\]
Integrating this equation gives
\[
\begin{array}{l}
h(x_0,e^{i\omega_1\xi}x_1,e^{-i\omega_1\xi}\overline{x_1},\ldots,e^{i\omega_p\xi}x_p,e^{-i\omega_p\xi}\overline{x_p},\nu,\mu)-h(x_0,x_1,\overline{x_1},\ldots,x_p,\overline{x_p},\nu,\mu)=\\
\\
\xi\,\mathfrak{g}(x_0,x_1,\overline{x_1},\ldots,x_p,\overline{x_p},\nu,\mu).
\end{array}
\]
For any given $(x,\nu,\mu)\in\mathbb{R}^{2p+2+s}$, the left hand side
of the previous equation is bounded in $\xi$, and it follows that
$\mathfrak{g}$ must equal $0$ and that
\[
h(x_0,e^{i\omega_1\xi}x_1,e^{-i\omega_1\xi}\overline{x_1},\ldots,e^{i\omega_p\xi}x_p,e^{-i\omega_p\xi}\overline{x_p},\nu,\mu)\equiv
h(x_0,x_1,\overline{x_1},\ldots,x_p,\overline{x_p},\nu,\mu).
\]
The conclusion follows from the non-resonance condition on
$\omega_1,\ldots,\omega_p$ specified in Hypothesis \ref{spectralhyp},
using density and continuity.
\hfill
\qed
\begin{lemma}
Equation (\ref{2ndcompnts}) implies that $h^{\nu}$ must be of the form
\[
\begin{array}{lll}
h^{\nu}=\mathfrak{h}(x_0,x_1\overline{x_1},\ldots,x_p\overline{x_p},\nu,\mu)&=&{\displaystyle
-\int_0^{x_0}\,r_1(t,x_1\overline{x_1},\ldots,x_p\overline{x_p},\mu)\,dt}\,-\\&&\\&&
\mathfrak{m}(x_1\overline{x_1},\ldots,x_p\overline{x_p},\nu,\mu).
\end{array}
\]
\label{applem2}
\end{lemma}
\proof
Write $h^{\nu}(x,\nu,\mu)=\sum_{j=0}^{\ell}\,a_j(x,\mu)\nu^j$.
Equation (\ref{2ndcompnts}) becomes
\begin{equation}
\sum_{j=0}^{\ell}\,\frac{\partial a_j}{\partial x_0}\nu^{j+1}+\sum_{j=0}^{\ell}\,\left(\sum_{j=1}^p\,i\omega_j
\left(x_j\frac{\partial a_j}{\partial
    {x_j}}-\overline{x_j}\frac{\partial a_j}{\partial
    {\overline{x_j}}}\right)\right)\nu^j=\mathfrak{g},
\label{wassuppp3}
\end{equation}
where $\mathfrak{g}=-(\nu\,r_1+r_2)$ is $\mathbb{T}^p$ invariant as in
Lemma \ref{applem1}.  Applying
Lemma \ref{applem1} successively to the coefficient of $\nu^0$,
$\nu^1,\ldots,\nu^{\ell+1}$ in (\ref{wassuppp3}), we get that $r_2=0$, 
${\dps\frac{\partial a_0}{\partial x_0}=-r_1}$, ${\dps\frac{\partial
    a_j}{\partial x_0}=0, j=1,\ldots,\ell}$, and $a_0,\ldots,a_{\ell}$
    are $\mathbb{T}^p$-invariant.  The conclusion of the lemma then follows
    immediately upon integration.
\hfill
\qed

The component of (\ref{1stcompnt}) along
$\mathbf{e}_0$ now has the form
\begin{equation}
\begin{array}{lll}
f_{0}&=&{\dps\nu\frac{\partial h^x_{0}}{\partial {x_0}}+
\sum_{j=1}^p\,i\omega_j\left(x_j\,\frac{\partial h^x_{0}}{\partial
    {x_j}}-\overline{x_j}\,\frac{\partial h^x_{0}}{\partial {\overline{x_j}}}\right)
+}\\&&\\
&&{\displaystyle
\int_0^{x_0}\,r_1(t,x_1\overline{x_1},\ldots,x_p\overline{x_p},\mu)\,dt}+
\mathfrak{m}(x_1\overline{x_1},\ldots,x_p\overline{x_p},\nu,\mu)+\\&&\\
&&x_0\,r_1(x_0,x_1\overline{x_1},\ldots,x_p\overline{x_p},\mu).
\end{array}
\label{1stcompnt2}
\end{equation}
Using Taylor's theorem, we write
\[
\mathfrak{m}(x_1\overline{x_1},\ldots,x_p\overline{x_p},\nu,\mu)=
\mathfrak{m}(x_1\overline{x_1},\ldots,x_p\overline{x_p},0,\mu)+\nu\,\hat{\mathfrak{m}}
(x_1\overline{x_1},\ldots,x_p\overline{x_p},\nu,\mu),
\]
and note that 
\[
\begin{array}{l}
{\dps\nu\,\hat{\mathfrak{m}}
(x_1\overline{x_1},\ldots,x_p\overline{x_p},\nu,\mu)=
\nu\frac{\partial}{\partial {x_0}}\,x_0\hat{\mathfrak{m}},\,\,\,\,\,\,\mbox{\rm and}}\\ \\{\displaystyle
\sum_{j=1}^p\,i\omega_j\left(x_j\frac{\partial}{\partial
    {x_j}}\,x_0\hat{\mathfrak{m}}-\overline{x_j}\frac{\partial}{\partial {\overline{x_j}}}\,x_0\hat{\mathfrak{m}}\right)=0.}
\end{array}
\]
Consequently, if we denote
$\widehat{h^x_{0}}=h^x_{0}+x_0\hat{\mathfrak{m}}$, then (\ref{1stcompnt2})
reduces to
\begin{equation}
\begin{array}{lll}
f_{0}&=&{\dps\nu\,\frac{\partial\widehat{h^x_{0}}}{\partial x_0}+
\sum_{j=1}^p\,i\omega_j\left(x_j\,\frac{\partial\widehat{h^x_{0}}}{\partial
    x_j}-\overline{x_j}\,\frac{\partial\widehat{h^x_{0}}}{\partial\overline{x_j}}\right)}
+\\&&\\
&&{\displaystyle
\int_0^{x_0}\,r_1(t,x_1\overline{x_1},\ldots,x_p\overline{x_p},\mu)\,dt}+
\mathfrak{m}(x_1\overline{x_1},\ldots,x_p\overline{x_p},0,\mu)+\\&&\\
&&x_0\,r_1(x_0,x_1\overline{x_1},\ldots,x_p\overline{x_p},\mu).
\end{array}
\label{1stcompnt3}
\end{equation}
Together, (\ref{1stcompnt3}) and the last $2p$ components of (\ref{1stcompnt})
imply that $f$ can be written as a sum of an element in $\mbox{\rm
  range}\,{\cal L}_{{B},\nu}$ and an element in
$H^{2p+1+s}_{\ell}(\mathbb{R}^{2p+1},\mathbb{T}^p)$, i.e.
\[
H^{2p+2+s}_{\ell}(\mathbb{R}^{2p+1})=H^{2p+1+s}_{\ell}(\mathbb{R}^{2p+1},\mathbb{T}^p)+
\mbox{\rm range}\,{\cal L}_{{B},{\nu}}.
\]

We now must show that the above sum is, in fact, a direct sum.
Suppose that $h^x\in H^{2p+2+s}_{\ell}(\mathbb{R}^{2p+1})$
and $g^x\in H^{2p+1+s}_{\ell}(\mathbb{R}^{2p+1},\mathbb{T}^p)$
are such that
\begin{equation}
{\cal L}_{{B},\nu}\,h^x\,+\,g^x=0.
\label{case2homeq}
\end{equation}
Write the component of $g^x$ along
$\mathbf{e}_0$ as
\begin{equation}
g^x_{0}(x_0,x_1\overline{x_1},\ldots,x_p\overline{x_p},\mu)=
g^x_{0}(0,x_1\overline{x_1},\ldots,x_p\overline{x_p},\mu)+
x_0\widehat{g^x_{0}}(x_0,x_1\overline{x_1},\ldots,x_p\overline{x_p},\mu),
\label{whatup1}
\end{equation}
and define
\[
\begin{array}{l}
r_1(x_0,x_1\overline{x_1},\ldots,x_p\overline{x_p},\mu)=\\\\{\displaystyle\frac{1}{x_0^2}\,\int_0^{x_0}\,
\left(t\widehat{g^x_{0}}(t,x_1\overline{x_1},\ldots,x_p\overline{x_p},\mu)+
t^2\,\frac{\partial\widehat{g^x_{0}}}{\partial {x_0}}
(t,x_1\overline{x_1},\ldots,x_p\overline{x_p},\mu)\right)\,dt.}
\end{array}
\]
Note that $r_1$ is regular at $x_0=0$, since both the integral term above
and its derivative with respect to $x_0$ vanish at $x_0=0$.  A
simple computation verifies that 
\begin{equation}
\begin{array}{lll}
x_0\,\widehat{g^x_{0}}(x_0,x_1\overline{x_1},\ldots,x_p\overline{x_p},\mu)&=&{\displaystyle
\int_0^{x_0}r_1(t,x_1\overline{x_1},\ldots,x_p\overline{x_p},\mu)\,dt+}\\&&\\&&
x_0\,r_1(x_0,x_1\overline{x_1},\ldots,x_p\overline{x_p},\mu).
\end{array}
\label{whatup2}
\end{equation}
Consequently, if we define $g^{\mu}=0$, $h^{\mu}=0$,
\[
g^{\nu}=\nu\,r_1(x_0,x_1\overline{x_1},\ldots,x_p\overline{x_p},\mu)
\]
and
\[
h^{\nu}=-\int_0^{x_0}\,r_1(t,x_1\overline{x_1},\ldots,x_p\overline{x_p},\mu)\,dt-
g^x_{0}(0,x_1\overline{x_1},\ldots,x_p\overline{x_p},\mu),
\]
then from Lemma \ref{lemenf1}, we see that $((x_0\,r_1,g^x_1,\ldots,g^x_{2p}),g^{\nu},g^{\mu})\in H^{2p+2+s}_{\ell}(\mathbb{R}^{2p+2+s},\Gamma)$.
It now follows from (\ref{case2homeq}) that
\[
(0,0,0)={\cal L}_{\tilde{B}}(h^x,h^{\nu},h^{\mu})+((x_0\,r_1,g^x_1,\ldots,g^x_{2p}),g^{\nu},g^{\mu}).
\]
From Proposition \ref{prop_enf1}, we get that $r_1=g^x_{1}=\cdots
=g^x_{2p}=0$ and ${\cal L}_{\tilde{B}}(h^x,h^{\nu},h^{\mu})=(0,0,0)$,
from which it follows that
\[
{\cal
  L}_{{B},\nu}\,h^x+g^x_{0}(0,x_1\overline{x_1},\ldots,x_p\overline{x_p},\mu)\,\mathbf{e}_{0}=0.
\]
The component of the above equation along
$\mathbf{e}_0$, is
\[
\nu\,\frac{\partial h^x_{0}}{\partial
  x_0}+\sum_{j=1}^p\,i\omega_j\left(x_j\,\frac{\partial
  h^x_{0}}{\partial x_j}-\overline{x_j}\,\frac{\partial
  h^x_{0}}{\partial\overline{x_j}}\right)=-
g^x_{0}(0,x_1\overline{x_1},\ldots,x_p\overline{x_p},\mu).
\]
Using Lemma \ref{applem1}, we conclude that
$g^x_{0}(0,x_1\overline{x_1},\ldots,x_p\overline{x_p},\mu)=0$,
${\cal L}_{{B},\nu}\,h^x=0$,
and from (\ref{whatup1}), that $g^x=0$.
Therefore, we conclude
that
\[
H^{2p+2+s}_{\ell}(\mathbb{R}^{2p+1})=H^{2p+1+s}_{\ell}(\mathbb{R}^{2p+1},\mathbb{T}^p)\oplus
\mbox{\rm range}\,{\cal L}_{{B},{\nu}}.
\]
\hfill
\qed

\Section{Proof of Proposition~\ref{prop_Adecomposition}}

The proof that $\mathring{A}$ is a projection is similar to the proof given for
the projection operator $A$ used in the semisimple case in \cite{CL1}.

Now, let $f\in\mbox{\rm range}\,\, \mathring{A}$, then $\mathring{A}f=f$, i.e.
\[
f(x,\nu,\mu)=
\int_{\Gamma_0}\,\gamma\,f(\gamma^{-1}x,0,\mu)\,d\gamma.
\]
So, for any $\sigma\in\Gamma_0$, we have
\[
\begin{array}{lll}
\sigma\,f(\sigma^{-1}x,\nu,\mu)&=&
{\displaystyle\sigma\,\int_{\Gamma_0}\,\gamma\,f(\gamma^{-1}\sigma^{-1}x,0,\mu)\,d\gamma
  =
  \int_{\Gamma_0}\,\sigma\gamma\,f((\sigma\gamma)^{-1}x,0,\mu)\,d\gamma}\\[0.15in]
&=&{\displaystyle\int_{\Gamma_0}\,\gamma\,f(\gamma^{-1}x,0,\mu)\,d\gamma =
  f(x,\nu,\mu)=f(x,0,\mu)},
\end{array}
\]
and thus $f\in H^{2p+1+s}_{\ell}(\mathbb{R}^{2p+1},\mathbb{T}^p)$.  
On the other hand, if $f\in
H^{2p+s+1}_{\ell}(\mathbb{R}^{2p+1},\mathbb{T}^p)$, then
\[
\begin{array}{l}
{\displaystyle(\mathring{A}f)(x,\nu,\mu)=\int_{\Gamma_0}\,\gamma\,f(\gamma^{-1}x,0,\mu)\,d\gamma=
\int_{\Gamma_0}\,f(x,0,\mu)\,d\gamma=}\\\\f(x,0,\mu)=f(x,\nu,\mu),
\end{array}
\]
so $f\in\mbox{\rm range}\,\, \mathring{A}$.  
This establishes (\ref{Arange}).
We now establish (\ref{Aker}).  Since $\mathring{A}$ is a projection, then
\[
H^{2p+2+s}_{\ell}(\mathbb{R}^{2p+1})=\mbox{\rm
  range}\,\mathring{A}\oplus\mbox{\rm ker}\,\mathring{A}.
\]
From Proposition \ref{prop_enf2}, we conclude that $\mbox{\rm dim
  ker}\,\mathring{A}=\mbox{\rm dim range}\,{\cal L}_{{B},\nu}$.  Thus, we need only show that
$\mbox{\rm range}\,{\cal L}_{{B},\nu}\subset\mbox{\rm
  ker}\,\mathring{A}$.  
Recall the following lemma which was proved in \cite{CL1}.
\begin{lemma}
Let $g:\Gamma_0\longrightarrow\mathbb{R}^{2p+1}$ be a continuous
function, then
\[
\int_{\Gamma_0}\,g(\gamma)\,d\gamma =
\lim_{T\rightarrow\infty}\,\frac{1}{T}\,\int_0^T\,g(e^{{B}s})\,ds.
\]
\label{averaging_lemma}
\end{lemma}

Now, let $f\in\mbox{\rm range}\,{\cal L}_{{B},\nu}$, then there exists $g\in
H^{2p+2+s}_{\ell}(\mathbb{R}^{2p+1})$ such that
\[
D_xg(x,\nu,\mu){B}x-{B}g(x,\nu,\mu)+\nu\,\frac{\partial
  g}{\partial {x_0}}(x,\nu,\mu)=f(x,\nu,\mu),\,\,\,\forall\,(x,\nu,\mu)\in\mathbb{R}^{2p+2+s}.
\]
Therefore, using Lemma \ref{averaging_lemma}, we get
\[
\begin{array}{l}
{\displaystyle
(\mathring{A}f)(x,\nu,\mu)=\int_{\Gamma_0}\,\gamma\,f(\gamma^{-1}x,0,\mu)\,d\gamma =
\lim_{T\rightarrow\infty}\,\frac{1}{T}\int_0^T\,e^{{B}s}\,f(e^{-{B}s}x,0,\mu)\,ds}\\[0.15in]
{\displaystyle
  =\lim_{T\rightarrow\infty}\,\frac{1}{T}\,\int_0^T\,e^{{B}s}\,\left(D_xg(e^{-{B}s}x,0,\mu)
{B}e^{-{B}s}x-{B}g(e^{-{B}s}x,0,\mu)\right)\,ds}\\[0.15in]
{\displaystyle
  =\lim_{T\rightarrow\infty}\,\frac{1}{T}\,\int_0^T\,\frac{d}{ds}
\left(-e^{{B}s}g(e^{-{B}s}x,0,\mu)\right)\,ds}\\[0.15in]
{\displaystyle
=\lim_{T\rightarrow\infty}\,\frac{-e^{{B}T}g(e^{-{B}T}x,0,\mu)+g(x,0,\mu)}{T}}
\end{array}
\]
and this last limit is equal to $0$, since the numerator is bounded
in $T$ for any given $(x,\mu)\in\mathbb{R}^{2p+1+s}$.
So we conclude that $f\in\mbox{\rm ker}\,A$, and thus that
$\mbox{\rm ker}\,A=\mbox{\rm range}\,{\cal L}_{{B},\nu}$.  This establishes (\ref{Aker}), and
concludes the proof of Proposition \ref{prop_Adecomposition}.
\hfill\qed


\begin{thebibliography}{10}
\bibitem{BC94} 
J.~B\'elair and S.A.~Campbell. 
\newblock Stability and bifurcations of equilibria in a multiple-delayed
differential equation.
\newblock {\em SIAM J.~Appl.~Math.} {\bf 54}, (1994) 1402--1424.

\bibitem{BBL} 
A.~Beuter, J.~B\'elair and C.~Labrie. 
\newblock Feedback and delays in neurological diseases : a modeling study 
using dynamical systems.
\newblock {\em Bulletin Math. Biology} {\bf 55}, (1993) 525--541.

\bibitem{BB} 
P-L.~Buono and J.~B\'elair.
\newblock Restrictions and unfolding of double Hopf bifurcation in 
functional differential equations. 
\newblock {\em J.~Diff.~Eqs.} {\bf 189}, (2003) 234--266.

\bibitem{CL1} 
Y.~Choi and V.G.~LeBlanc. 
\newblock Toroidal normal forms for bifurcations in retarded
functional differential equations I: Multiple Hopf and
transcritical/multiple Hopf interaction. 
\newblock {\em Submitted} (2005), preprint available at
arXiv:math.DS/0505392 v1.

\bibitem{ETBCI}
C.~Elphick, E.~Tirapegui, M.E.~Brachet, P.~Coullet and G.~Iooss.
\newblock A simple global characterization for normal forms of
singular vector fields.
\newblock {\em Phys. D.} {\bf 29}, (1987) 95--127.

\bibitem{FMTB}
T. Faria and L.T. Magalh$\tilde{\mbox{\rm a}}$es. 
\newblock Normal Forms for Retarded Functional
Differential Equations and Applications to Bogdanov-Takens Singularity.
\newblock {\em J. Diff. Eqs.} {\bf 122}, (1995) 201--224.

\bibitem{FMH}
T. Faria and L.T. Magalh$\tilde{\mbox{\rm a}}$es. 
\newblock Normal Forms for Retarded Functional
Differential Equations with Parameters and Applications to Hopf Bifurcation.
\newblock {\em J. Diff. Eqs.} {\bf 122}, (1995) 181--200.

\bibitem{FMR}
T. Faria and L.T. Magalh$\tilde{\mbox{\rm a}}$es. 
\newblock Realisation of Ordinary Differential
Equations by Retarded Functional Differential Equations in Neighborhoods of
Equilibrium Points.
\newblock {\em Proc. Roy. Soc. Ed.} {\bf 125A}, (1995) 759--776.

\bibitem{FM96} 
T. Faria and L.T. Magalh$\tilde{\mbox{\rm a}}$es.
\newblock Restrictions on the possible flows of scalar retarded
functional differential equations in neighborhoods of singularities.
\newblock {\em J.~Dyn.~Diff.~Eqs} {\bf 8}, (1996) 35--70.

\bibitem{GSSI}
M.~Golubitsky and D.G.~Schaeffer. 
\newblock {\em Singularities and Groups in
  Bifurcation Theory. Vol. 1.}, Applied Mathematical Sciences {\bf
  51}, Springer-Verlag, New York, (1985).

\bibitem{GSSII}
M.~Golubitsky, I.~Stewart and D.G.~Schaeffer. 
\newblock {\em Singularities and Groups in Bifurcation Theory,
  Vol. 2.}, Applied Mathematical Sciences {\bf 69}, Springer-Verlag,
New York, (1988).

\bibitem{Guck1}
J.~Guckenheimer.
\newblock On a codimension two bifurcation.
\newblock {\em Lecture Notes in Mathematics} {\bf 898}, (1981)
99--142.

\bibitem{Guck2}
J.~Guckenheimer.
\newblock Multiple bifurcations of codimension two.
\newblock {\em SIAM J. Math. Anal.} {\bf 15}, (1984) 1--49.

\bibitem{Hal85}
J.K.~Hale.
\newblock Flows on center manifolds for scalar functional differential
equations.
\newblock {\em Proc. Roy. Soc. Edinburgh} {\bf 101}, (1985) 193--201.

\bibitem{HalVL}
J.K. Hale and S.M. Verduyn Lunel. 
\newblock {\em Introduction to Functional Differential
Equations}.
\newblock Springer-Verlag, New York, (1993).

\bibitem{HFEKGG}
T.~Heil, I.~Fischer, W.~Elsa{\ss}er, B.~Krauskopf, K.~Green and A.~Gavrielides.
\newblock Delay dynamics of semiconductor lasers with short external
cavities: Bifurcation scenarios and mechanisms. 
\newblock {\em Phys. Rev. E} {\bf 67}, (2003) 066214-1--066214-11.

\bibitem{Kuang}
Y.~Kuang.
\newblock {\em Delay differential equations with applications in population 
dynamics.} 
\newblock Mathematics in Science and Engineering, 191. Academic Press,  
Boston, (1993).

\bibitem{LM} 
A.~Longtin and J.G.~Milton. 
\newblock Modelling autonomous oscillations in the human pupil light 
reflex using nonlinear delay-differential equations.
\newblock {\em Bulletin Math. Biology} {\bf 51}, (1989) 605--624.

\bibitem{SC} 
E.~Stone and S.A.~Campbell.
\newblock Stability and bifurcation analysis of a nonlinear DDE model
for drilling. 
\newblock {\em J. Nonlinear Sci.} {\bf 14} (2004), 27--57.

\bibitem{SS} 
M.J. Suarez and P.L. Schopf. 
\newblock A Delayed Action Oscillator for ENSO. 
\newblock
{\em J. Atmos. Sci.}  {\bf 45} (1988), 3283--3287.


\end{thebibliography}
\end{document}